\documentclass[11pt]{article}
\usepackage{amssymb}
\usepackage{amsmath}
\usepackage{latexsym}
\usepackage{hyperref}

\usepackage[all]{xy}
\usepackage{color}
\usepackage[normalem]{ulem}

\numberwithin{equation}{section}
\newtheorem{thm}[equation]{Theorem}

\newtheorem{defn}[equation]{Definition}
\newtheorem{rem}[equation]{Remark}

\newtheorem{corol}[equation]{Corollary}

\title{On an interior Calder\'{o}n  operator and  a related Steklov eigenproblem for Maxwell's equations\footnote{Submitted for publication to  Siam Journal on Mathematical Analysis on 21 March 2019, revised on 12 May 2020, accepted for publication  on 16 July 2020.  }}
\author{Pier Domenico Lamberti and Ioannis G. Stratis}

\date{\today }

\begin{document}

\newcommand{\N}{\mathbb{N}}
\newcommand{\R}{\mathbb{R}}


\newcommand{\cu}{{\rm {curl}}}
\newcommand{\di}{{\rm div}}
\newcommand{\gr}{{\rm grad}}
\newcommand{\dig}{{\rm div}_{\scriptscriptstyle\Gamma}}
\newcommand{\grg}{{\rm grad}_{\scriptscriptstyle\Gamma}}
\newcommand{\cug}{{\rm curl}_{\scriptscriptstyle\Gamma}}
\newcommand{\ca}{ {\mathcal{C}}_{\rm me}}
\newcommand{\cas}{{\mathcal{C}}_{\rm \tiny H}^{{\mathcal{\tiny S}}}}
\newcommand{\cald}{{\mathcal{C}}_{\rm \tiny EH}}
\newcommand{\caH}{{\mathcal{C}}_{\rm \tiny H}}
\newcommand{\cau}{{\mathcal{C}}}

\maketitle
\begin{center}
\texttt{Dedicated to Professor Nicholas D. Alikakos\\ on the occasion of his retirement  }
\end{center}
\vspace{0.4cm}

\noindent
{\bf Abstract:} We discuss a Steklov-type problem for Maxwell's equations  which is related to an interior Calder\'{o}n operator and an appropriate  Dirichlet-to-Neumann  type map. The corresponding Neumann-to-Dirichlet map turns out to be compact and this provides a Fourier basis of Steklov eigenfunctions for the associated energy spaces.  With an   approach  similar to that developed by  Auchmuty for the Laplace operator,  we provide natural spectral representations for the appropriate trace spaces, for the Calder\'{o}n operator itself and
for the solutions of the corresponding boundary value problems subject to electric or magnetic boundary conditions on a cavity.

\vspace{11pt}


\vspace{6pt}

\section{Introduction}

In a homogeneous isotropic medium filling a domain $\Omega$ in $\R^3$ the {\it time-harmonic} Max\-well's equations read
\begin{equation}\label{Maxwelltimeharm}
\cu E - {\rm i} \, \omega \, \mu \, H = 0\,,\,\,\cu H + {\rm i} \, \omega \,\varepsilon \, E = 0\,,
\end{equation}
where $E, H$ are, respectively, the spatial parts of the electric and the magnetic field, $\varepsilon$ and $\mu$ are the electric permittivity and the magnetic permeability of the medium, and  $\omega  > 0$ is the angular frequency (we have adopted the time convention $\mathrm{e}^{- \,\mathrm{i}\, \omega\, t}$). In the considered case of homogeneous isotropic media, $\varepsilon$ and $\mu$ are constants, therefore $E$ and $H$ are automatically divergence-free.

Of {\it sine qua non} importance in electromagnetics is the following boundary value problem, involving the so-called ``perfect conductor'' condition on the boundary $\Gamma$ of $\Omega$, i.e.,  the tangential trace of the electric field on  $\Gamma$ is given by a fixed vector $m$:
\begin{equation}\label{PEC}
\left\{
\begin{array}{ll}
\cu E - {\rm i} \, \omega \, \mu \, H = 0\,,\,\, \cu H + {\rm i} \, \omega \, \varepsilon \, E = 0 ,& \ \ {\rm in}\ \Omega ,\\
\nu \times E = m ,& \ \ {\rm on}\ \Gamma .
\end{array}
\right.
\end{equation}

The interior Calder\'on operator is defined as the mapping of the tangential component of the electric field to the tangential component of the magnetic field on $\Gamma$, i.e., $m \mapsto \nu \times H$. This is the origin of another term used for this operator, namely the ``electric to magnetic boundary component map''. Calder\'on operators are also known as capacity, or impedance, or admittance, or
Poincar\'{e}-Steklov operators.

 Operating by $\cu$ on \eqref{PEC} and setting $\widetilde{m} := - \,{\rm i} \, \omega \,\varepsilon \, m$ we obtain
\begin{equation}\label{PECH}
\left\{
\begin{array}{ll}
\cu \cu H - \, \omega^2  \varepsilon \, \mu \, H = 0 ,& \ \ {\rm in}\ \Omega ,\\
\nu \times \cu H =  \widetilde{m},& \ \ {\rm on}\ \Gamma .
\end{array}
\right.
\end{equation}
The corresponding interior Calder\'on operator for \eqref{PECH}, maps $\widetilde{m}$ to $ \nu \times H$.

Let us note that, in view of the vector identity
$ \cu \, \cu w =  \gr\,  \di w - \Delta w $ and  the fact that $H$ is divergence-free, \eqref{PECH} can also be written as
\begin{equation}\nonumber
\left\{
\begin{array}{ll}
\Delta H +  \omega^2 \, \varepsilon \, \mu \, H = 0 ,& \ \ {\rm in}\ \Omega ,\\
\nu \times \cu H = \widetilde{m} ,& \ \ {\rm on}\ \Gamma .
\end{array}
\right.
\end{equation}
Based on the sign of $\varepsilon$ and of $\mu$, materials can be classified as
\begin{itemize}
\item[(RH)] When $\varepsilon > 0,\, \mu > 0$ the materials are called ``right handed'' or ``double positive'', and exhibit forward propagating waves. These are the conventional materials of electromagnetics, e.g., dielectrics.
\item[(LH)] When $\varepsilon < 0,\, \mu < 0$ the materials are called ``left handed'' or ``double negative'', and they exhibit backward propagating waves. They are not found in Nature, but are physically realizable.
\item[(SRR)] The materials with $\varepsilon > 0,\, \mu < 0$ are called of ``split ring resonator'' structure, and exhibit evanescent decaying waves and no transmission (typical examples are ferrites, microstructured magnets and split rings).
\item[(TW)] The materials with $\varepsilon < 0,\, \mu > 0$ are called of ``thin wire'' structure, and they exhibit evanescent decaying waves and no transmission (typical examples are plasmas and fine wire structures).
\end{itemize}
(LH), (SRR) and (TW) are included in the so-called ``metamaterials'', \cite{li}.

\vspace{0.3cm}
Note that the results in the present work cover all the above cases (the coefficient $\alpha$ of our approach, that corresponds to the physical constant $\omega^2 \, \varepsilon \, \mu$ appearing in the time-harmonic Maxwell's equations, can be of any sign).

\vspace{0.3cm}

Let us note that the study of the operators ``$\,\cu\,\cu\,$'' and ``$\,\cu\,\cu - \varrho^2 I\,$'' is essential not only in the mathematical theory of classical electromagnetics, but in other related important applications areas, such as, e.g., the theory of superconductors, magnetohydrodynamics (MHD)(where the equations - consisting of an elegant and subtle coupling of the Navier-Stokes and Maxwell's equations - govern the motion of electrically conducting viscous incompressible fluids in a magnetic field (see, e.g., \cite{boumaam}), and in particular in the ideal linear MHD equations that describe the stability properties of a ``tokamak'', i.e., a device which uses a powerful magnetic field to confine a hot plasma in the shape of a torus), etc. 
\\ \\
The classical Steklov eigenproblem for a bounded smooth domain $\Omega$ in $\R^n$ with boundary $\Gamma$ reads
$$
\left\{
\begin{array}{ll}
\Delta v= 0,& \ \ {\rm in}\ \Omega,\\
D_{\nu}v =\lambda v ,& \ \ {\rm on}\ \Gamma,
\end{array}
   \right.
$$
where the  unknown $v$ is a real or complex-valued function called Steklov eigenfunction and the unknown $\lambda$ is a non-negative real number called Steklov eigenvalue.
Here $\nu$ denotes the unit outer normal to $\Gamma$ and $D_{\nu}v$ the  normal derivative of $v$.
The Steklov eigenvalues  can be equivalently defined as the eigenvalues of the celebrated DtN (Dirichlet-to-Neumann) map defined from $H^{1/2}(\Gamma)$ to
$H^{-1/2}(\Gamma)$ by $ g\mapsto  D_{\nu}v$ where $v$ is the solution to the Dirichlet problem
$$
\left\{
\begin{array}{ll}
\Delta v= 0,& \ \ {\rm in}\ \Omega,\\
v =g, & \ \ {\rm on}\ \Gamma .
\end{array}
   \right.
$$
It turns out that the non-zero eigenvalues are the reciprocals of the eigenvalues of the corresponding NtD (Neumann-to-Dirichlet) map which can be considered as a compact self-adjoint map from
$L^2(\Gamma)$ to itself. In particular, the eigenvalues have finite multiplicity and can be represented as a non-decreasing divergent sequence.
We refer  to \cite{girpol,  lampro1}  for an introduction to  Steklov-type problems,  and to \cite{hey} for  an interesting  application of the problem.

Although the study of Steklov eigenvalues has a long history (see \cite{kukukwnapoposi}) and  Steklov boundary conditions  have been considered for many classes of operators, in the literature there are not so many results concerning Steklov-type eigenvalues for Maxwell's equations; an interesting exception are those contained in the very recent papers \cite{calamo, cocomo}.

We note that in \cite{calamo} (along the approach introduced in \cite{cacomemo} for the Helmholtz equation) the use of Steklov eigenvalues for Maxwell's equations is suggested to detect changes in a scatterer using remote measurements of the scattered wave, i.e., as a novel ``target signature'' for nondestructive testing via inverse scattering. Because the Steklov eigenvalue problem for Maxwell's equations is not a standard eigenvalue problem for a compact operator, a modified Steklov problem is proposed, that restores compactness. In particular it is shown that it is possible to measure Steklov eigenvalues for a bounded inhomogeneous scatterer by solving a sequence of modified far field equations. To this end, the authors perturb the usual far field equation of the linear sampling method by using the far field pattern of an auxiliary impedance problem related to the modified Steklov problem. In order to measure the modified Steklov eigenvalues of a domain from far field measurements, the authors prove (i) the existence of modified Steklov eigenvalues, (ii) the well-posedness of the corresponding auxiliary exterior impedance problem, and (iii) provide theorems on the detection of modified Steklov eigenvalues from far field measurements. 
While our present paper was at the stage of review, two interesting manuscripts appeared in arXiv, continuing in a sense the work in \cite{calamo}: in \cite{halla1}, both the original Steklov eigenvalue problem for Maxwell's equations and  the aforementioned modified Steklov problem are studied, and their Fredholmness and approximation is analyzed. The original eigenvalue problem in the selfadjoint case is studied in \cite{halla2}, where it is established that, apart for a countable set of particular frequencies, the spectrum consists of three disjoint parts: the essential spectrum consisting of the point zero, an infinite sequence of positive eigenvalues which accumulate only at infinity, and an infinite sequence of negative eigenvalues which accumulate only at zero.  See also the very receent work by Cogar \cite{cogar-siam-2020jam}, \cite{cogar-arXiv_Jun-2020} and Cogar and Monk \cite{cogar-monk-arXiv_May-2020}.

The aim of the present paper is to furnish a natural Steklov problem for a class of Maxwell's equations, which intrinsically exhibits a discrete spectrum, and use it in the spirit of \cite{Auch2006} to provide spectral representations for the associated trace spaces and the solutions of the corresponding boundary value problems.
Our starting point is the observation that, in the mathematical theory of electromagnetism, a natural counterpart of the NtD map above is the aforementioned celebrated  Calder\'{o}n operator.
The Calder\'{o}n operator and its variants have been and are studied extensively, see, e.g., \cite{ces}. Here we focus on a Calder\'{o}n operator associated with the following boundary value problem: for a bounded, connected,   open set  $\Omega$ in $\R^3$ of class $C^{1,1}$ and  $\alpha \in \R$,  $\theta >0$ consider the interior\footnote{The term ``interior'' is used to emphasize  that the problem is cast in $\Omega$ and not in the exterior of $\Omega$, as is done in the case of  Calder\'{o}n operators for scattering problems.} problem

\begin{equation}\label{classiccucu1}
\left\{
\begin{array}{ll}
\cu\, \cu u -\alpha u -\theta\, \gr\,  \di u=0,& \ \ {\rm in}\ \Omega ,\\
\nu \times \cu u = f ,& \ \ {\rm on}\ \Gamma ,
\end{array}
\right.
\end{equation}
where $u$ is the unknown vector field.  
We point out  that the penalty term $\theta\, \gr\,  \di u$ is introduced in the equation in order to guarantee the coercivity of the quadratic form associated with the operator, as done e.g., in \cite{coda}. 
Note that the boundary operator  $\nu \times \cu u$ in \eqref{classiccucu1} can be considered  as the ``electromagnetic version'' of the boundary operator $D_{\nu}u$ on
$\Gamma$ usually associated with the scalar Laplace operator. Indeed, the boundary conditions in  \eqref{classiccucu1}  are the natural boundary conditions arising from the integration by parts formula
\begin{equation}\label{parts}
\int_{\Omega}\cu u\cdot  \cu \varphi\, dx =\int_{\Omega}\cu\, \cu u \cdot \varphi \,dx-\int_{\Gamma }(\nu \times \cu u )\cdot \varphi \, d\sigma\, ,
\end{equation}
which is valid for sufficiently regular vector fields $u, \varphi$, while the Neumann boundary conditions for the scalar  Laplace operator are the natural boundary conditions arising from the integration by parts formula  
$$
\int_{\Omega}\gr v\cdot  \gr \psi\, dx = -\int_{\Omega}\Delta v\, \psi \,dx   +  \int_{\Gamma }D_{\nu}v \,    \psi \, d\sigma\,,
$$
which, again, is valid for sufficiently  regular scalar functions $v, \psi$.

The interior Calder\'{o}n operator  $\cau$ is here defined by
\begin{equation}
\cau (f) = \nu \times u \, ,
\end{equation}
where $u$ is the solution of \eqref{classiccucu1}, see Section~\ref{calderonsec} for more details. Thus, the Calder\'{o}n operator establishes a correspondence between the electric and magnetic fields on  $\Gamma$ as follows:
$$
\nu\times \cu\, u\,  \longmapsto\,  \nu\times  u\, .
$$
Now, the corresponding to the NtD-map in the case of our problem  is 
$$
\nu\times \cu u\,  \longmapsto\,  u\,,
$$
the eigenvalues of which are the reciprocals of the eigenvalues of the following Steklov-type problem for Maxwell's equations
\begin{equation}\label{classiccucu1e}
\left\{
\begin{array}{ll}
\cu\, \cu u -\alpha u -\theta\, \gr\,  \di u=0,& \ \ {\rm in}\ \Omega ,\\
\nu \times \cu u = \lambda u ,& \ \ {\rm on}\ \Gamma .
\end{array}
\right.
\end{equation}
 It is important to clarify that the energy space involved  here is the classical space  $X_{\rm \scriptscriptstyle\tiny T}(\Omega) = H(\cu, \Omega) \cap H_{0}(\di, \Omega)$ and that problem
\eqref{classiccucu1e} has to be interpreted in the weak sense as follows\footnote{Note that the  minus sign in the right-hand side of \eqref{weakcucu1e} is due to the fact that the boundary condition in the problem is written in the customary form $\nu \times \cu u=\lambda u$, and it is responsible for the appearance of negative spectra.  The reader who is more familiar with classical eigenvalue problems for elliptic equations could write it  in the form $ \cu u\times \nu=\lambda u$, which  would  change the sign of the eigenvalues.
}:
find $u\in X_{\rm \scriptscriptstyle\tiny T}(\Omega)$  such that
\begin{equation}\label{weakcucu1e}
\int_{\Omega}\cu u\cdot  \cu \varphi\, dx -\alpha \int_{\Omega}u \cdot \varphi \, dx + \theta \int_{\Omega }\di u\,  \di \varphi \, dx =  -  \lambda \, \int_{\Gamma} u \cdot \varphi \, d\sigma\, ,
\end{equation}
for all $\varphi \in X_{\rm \scriptscriptstyle\tiny T}(\Omega)$. In particular,  any vector field $u\in X_{\rm \scriptscriptstyle\tiny T}(\Omega)$ satisfies the condition $u\cdot \nu =0 $ on $\Gamma$. 

 Note that the weak formulation \eqref{weakcucu1e} can be obtained from \eqref{classiccucu1e}  by a standard procedure as follows: assume that $u$ is a sufficiently regular  solution of problem \eqref{classiccucu1e}, say  $u\in (H^2(\Omega))^3$,    then multiplying both sides of the first equation in \eqref{classiccucu1e} by $\varphi \in X_{\rm \scriptscriptstyle\tiny T}(\Omega)$, integrating by parts and using formula \eqref{parts} and the boundary condition in  \eqref{classiccucu1e}, yield the validity of \eqref{weakcucu1e}.  By using the same calculations and the Fundamental Lemma of the Calculus of Variations, one can see that if   $u \in (H^2(\Omega))^3$ is     a  solution of $\eqref{weakcucu1e}$ then it is also a solution of $\eqref{classiccucu1e}$.  We also note that using the weak formulations allows to avoid assuming extra smoothness assumptions on the boundary of $\Omega$  required to guarantee the regularity of the solutions, see  e.g., \cite{weber}.

By the classical Gaffney inequality (Theorem~\ref{gaffney}), if $\Omega$ is of class $C^{1,1}$ the space  $X_{\rm \scriptscriptstyle\tiny T}(\Omega)$ is continuously embedded into $(H^1(\Omega))^3$.  Moreover, the trace operator is compact from $H^1(\Omega)$ to $L^2(\Gamma )$. It follows that if the
operator 
\[
u \mapsto \cu\, \cu u +\alpha u -\theta\, \gr\,  \di u
\] 
is invertible,
problem
\eqref{classiccucu1e} has a discrete spectrum and provides a Fourier basis of eigenfunctions for the space $TL^2(\Gamma)$ of square-summable vector fields tangent to $\Gamma$, see Section~\ref{eigensubsec}. For the sake of simplicity, these results are proved under the assumption that $\alpha <A_1 $ where $A_1>0$ is the first eigenvalue of the associated Dirichlet problem, see \eqref{ray} and \eqref{dirichlet}.  The analysis of the general case  is  discussed in  Section~\ref{app}, where we explain how to rule out the Dirichlet eigenfunctions associated with the Dirichlet eigenvalues smaller than $\alpha$.

We note that  the boundary condition appearing in the Steklov problem discussed in \cite[Equation~(2.6)]{calamo} is $\nu \times \cu u = \lambda u_{\rm \scriptscriptstyle\tiny T}$ on $\Gamma$, where $u_{\rm \scriptscriptstyle\tiny T}$ is the tangential component of $u$; accordingly the energy space used in \cite{calamo} to treat  that problem is    $ \{u\in H(\cu, \Omega):\ u_{\rm \scriptscriptstyle\tiny T}\in (L^2(\Gamma) )^3 \}$ and the corresponding eigenvectors turn out to be divergence free.
Our  boundary condition $\nu \times \cu u = \lambda u $ on $\Gamma$,  is clearly  stronger in the sense that it implies that our eigenfunctions are automatically tangential: this allows us to discard  the part of the spectrum associated with possible non tangential eigenvectors which  are responsible for the appearance of  an accumulation point in the spectrum of the operator discussed   in the counterexample  in  \cite[p.~4383]{calamo}. 

    With an approach analogous to that  in \cite{Auch2006} concerning the Laplace operator, our Steklov problem  allows us to achieve the following results:
 \begin{itemize}
 \item[(1)] we provide (Theorems~\ref{representation},~\ref{dualthm}) a spectral representation for the solutions of problem \eqref{classiccucu1},
 \item[(2)] we provide (see Theorem~\ref{representation} and Remark~\ref{finalrem}) a spectral representation for the  interior Calder\'{o}n operator $\cau$,
 \item[(3)] we provide (see Theorem~\ref{repdir}) a spectral representation for the  trace space $TH^{1/2}(\Gamma)$ and its  dual   $TH^{-1/2}(\Gamma)$,
 \item[(4)] we provide (see Theorem~\ref{repdir}) a spectral representation for the  solutions to the following  problem
\begin{equation}\label{weakcucudir}
\left\{
\begin{array}{ll}
\cu\, \cu u -\alpha u -\theta\, \gr\,  \di u=0,& \ \ {\rm in}\ \Omega ,\\
\nu \times  u = f ,& \ \ {\rm on}\ \Gamma ,\\
\nu \cdot u = 0 ,& \ \ {\rm on}\ \Gamma .
\end{array}
\right.
\end{equation}
 \end{itemize}

Another approach to the representation of an exterior Calder\'{o}n operator associated with a scattering problem for not necessarily spherical domains is proposed in \cite{kswy}; there the appropriate series expansions are performed with respect to generalized  harmonics (the set of eigenfunctions to the Laplace-Beltrami operator for the domain's boundary). Further, the norm in an appropriate trace space of the exterior Calder\'on operator is obtained in view of an eigenproblem for a suitable quadratic form. 
Let us note that in the case of a sphere the eigenfunctions of the Laplace-Beltrami operator are the spherical harmonics hence the classical Steklov eigenfunctions.
\vspace{0.3cm}

The present paper is organized as follows: Section~\ref{preliminaries} is devoted to preliminaries and notation; Section~\ref{interior} is devoted to the study of problems \eqref{classiccucu1} and \eqref{classiccucu1e}; in particular,  in Section~\ref{calderonsec} we introduce the Calder\'{o}n operator and the associated NtD-type  map. Section~\ref{eigensec} is  devoted to the  above mentioned spectral representations. Finally, Section~\ref{app} includes an appendix devoted to the study of the case $\alpha>A_1$.

\section{Preliminaries and notation}
\label{preliminaries}

Let $\Omega$ be a bounded open set  in $\R^3$ with  sufficiently smooth boundary $\Gamma :=\partial \Omega$. By $L^2(\Omega)$, $H^1(\Omega)$, $H^{1}_{0}(\Omega)$, $L^2(\Gamma)$, $H^{1/2}(\Gamma), H^{-1/2}(\Gamma)$,
we denote the standard Lebesgue and Sobolev spaces. We will also employ the following spaces:

\begin{itemize}
\item $H(\cu, \Omega) = \{u \in (L^2(\Omega))^3 : \cu u \in (L^2(\Omega))^3 \}\,,$\\
with norm: $ ||u||_{H(\cu, \Omega)} = \left( ||u||^{2}_{(L^2(\Omega))^3} + ||\cu u||^{2}_{(L^2(\Omega))^3} \right)^{1/2} $
\item $H(\di, \Omega) = \{u \in (L^2(\Omega))^3 : \di u \in L^2(\Omega) \}\,,$\\
with norm: $ ||u||_{H(\di, \Omega)} = \left( ||u||^{2}_{(L^2(\Omega))^3} + ||\di u||^{2}_{L^2(\Omega)} \right)^{1/2} $
\item $H_{0}(\di, \Omega) = \{u \in H(\di, \Omega): \nu \cdot u = 0 \,\, {\rm on}\ \Gamma \}$
\item $X_{\rm \scriptscriptstyle\tiny T}(\Omega) = H(\cu, \Omega) \cap H_0(\di, \Omega)\,,$\\
with norm:\! $ ||u||_{H(\cu, \Omega) \cap H(\di, \Omega)}\! =\! \left( ||u||^{2}_{(L^2(\Omega))^3} + ||\cu u||^{2}_{(L^2(\Omega))^3} + ||\di u||^{2}_{L^2(\Omega)} \right)^{1/2} $
\item $X_{\rm \scriptscriptstyle\tiny T}(\di \, 0, \Omega) = \{u \in X_{\rm \scriptscriptstyle\tiny T}(\Omega) : \di u = 0 \,\, {\rm in}\ \Omega \}$
\item $TL^2(\Gamma) =  \{u \in (L^2(\Gamma))^3 : \nu \cdot u = 0 \,\, {\rm on}\ \Gamma \}$
\item $TH^{1/2}(\Gamma) =  \{u \in (H^{1/2}(\Gamma))^3 : \nu \cdot u = 0 \,\, {\rm on}\ \Gamma \}$
\item $TH^{-1/2}(\Gamma) =(  TH^{1/2}(\Gamma)    )'     $
\end{itemize}

By $\di_{\scriptscriptstyle\Gamma}$ and $\gr_{\scriptscriptstyle\Gamma}$ we denote the usual tangential operators.\\
For details on  these operators and  spaces we refer to \cite{acl}, \cite{ces}, \cite{dali3}, \cite{gira}, \cite{kihe}, \cite{monk}, \cite{rsy}.
\newline

 Throughout this paper, we  consider bounded open sets of class $C^{1,1}$ in which case   $X_{\rm \scriptscriptstyle\tiny T}(\Omega)$ is continuously embedded in $(H^1(\Omega))^3$, and compactly embedded in $(L^2(\Omega))^3$. In particular, the following theorem holds (cf. e.g. \cite[Theorem~3.8]{gira}).

\begin{thm}\label{gaffney} If $\Omega $ is a bounded open set in $\R^3$ of class $C^{1,1}$
 then $X_{\rm \scriptscriptstyle\tiny T}( \Omega)$ is continuously embedded in $(H^1(\Omega ))^3$, in particular there exists
$c>0$ such that
$$\| u\|_{(H^1(\Omega))^3}\le c \, \left(       \|   u\|_{L^2(\Omega)}  + \|  \cu u\|_{L^2(\Omega)} +  \|  \di u\|_{L^2(\Omega)}  \right),$$
for all $u\in   X_{\rm \scriptscriptstyle\tiny T}( \Omega)$.
\end{thm}

Thus, for bounded  open sets of class $C^{1,1}$ we have that\footnote{   Although this space could be denoted with other more specific symbols, such as $TH^1(\Omega)$, we prefer to keep the above  notation  }
$$
X_{\rm \scriptscriptstyle\tiny T}( \Omega)=\left\{u\in  (H^1(\Omega))^3:\ u\cdot \nu =0\right\}\, .
$$

For a smooth vector field $w$ defined on $\overline{\Omega}$, its {\it ``tangential trace''} $\nu \times w$ on $\Gamma$ is denoted by $\gamma_{\rm \scriptscriptstyle\tiny T}w$, while its {\it ``tangential components trace''} $\nu \times (w \times \nu)$ on $\Gamma$ is denoted by $\pi_{\rm \scriptscriptstyle\tiny T}w$.

Since we are primarily interested in $\pi_{\rm \scriptscriptstyle\tiny T}$, we note       that $\pi_{\rm \scriptscriptstyle\tiny T}\left( (H^1(\Omega))^3 \right) = TH^{1/2}(\Gamma)$, $\pi_{\rm \scriptscriptstyle\tiny T}\left( (H^{1/2}(\Gamma))^3 \right) = TH^{1/2 }(\Gamma)$ and $\left( \pi_{\rm  \scriptscriptstyle\tiny T}\left( (H^{1/2}(\Gamma))^3 \right) \right)' = TH^{-1/2}(\Gamma)$.

Moreover,  we have that $\pi_{\rm \scriptscriptstyle\tiny T}$ is a compact operator from $(H^1(\Omega))^3$ to $TL^2(\Gamma)$. Note   that we shall often use the same symbol for a function and its  trace.

In Section~\ref{eigensec} we shall give an equivalent definition of the spaces $TH^{1/2}(\Gamma)$ and $TH^{-1/2}(\Gamma)$ based on an intrinsic spectral representation.

\section{Interior problems}
\label{interior}

Let $\Omega$ be a bounded, connected open set in $\R^3$ with boundary $\Gamma :=\partial \Omega \in C^{1,1}$. Let $\alpha \in \mathbb{R}$ and $\theta >0$ be fixed.
In this section, we discuss the solvability of problems \eqref{classiccucu1} and \eqref{classiccucu1e}.  We begin by proving that
 for any fixed  $\eta \geq 0$ sufficiently large, the problem
\begin{equation}\label{classiccucu1conc}
\left\{
\begin{array}{ll}
\cu\, \cu u -\alpha \ u -\theta\, \gr\,  \di u=0,& \ \ {\rm in}\ \Omega ,\\
\nu \times \cu u -\eta \, u = f ,& \ \ {\rm on}\ \Gamma
\end{array}
\right.
\end{equation}
has a (unique) solution for every datum $f \in TL^2(\Gamma)$.

The weak formulation  of problem \eqref{classiccucu1conc} reads: find  $u\in X_{\rm \scriptscriptstyle\tiny T}(\Omega )$ such that
\begin{equation}\label{weakcucu1}
\int_{\Omega}\cu u \cdot \cu \varphi\, dx \, -\, \alpha  \int_{\Omega}u\cdot \varphi\, dx \, + \, \theta  \int_{\Omega }\di u\,  \di \varphi \, dx \, +\, \eta \int_{\Gamma} u\cdot \varphi \, d\sigma =  - \int_{\Gamma} f \cdot\varphi \, d\sigma\, ,
\end{equation}
for all $\varphi \in X_{\rm \scriptscriptstyle\tiny T}(\Omega )$. Note that the equivalence of the two formulations for smooth solutions  can be proved in the same way explained in the Introduction for the case $\eta =0$. (Recall that using
 the weak formulation enables us to avoid extra regularity assumptions on $\Omega$.) 

 We note that when $\alpha >0$ the quadratic form associated with the left-hand side of equation  \eqref{weakcucu1}
is not necessarily positive and this complicates the analysis of the problem. In order to avoid a number of technicalities which would 
render the exposition much heavier, we find it not only simpler, but also natural to assume that $\alpha <A_1$, where 
\begin{equation}
\label{ray}
A_1:=\inf_{\substack{u\in (H^1_0(\Omega))^3\\u\ne 0 }  }\frac{  \int_{\Omega}|\cu u|^2\, dx + \theta \, \int_{\Omega }|\di u|^2 \, dx   }{\int_{\Omega} |u|^2dx}. 
\end{equation}
See Section~\ref{app} for a more general condition allowing arbitrarily large values of $\alpha$.

Note that $A_1>0$ by Theorem~\ref{gaffney},  and that $A_1$ is the first eigenvalue of the problem
\begin{equation}\label{dirichlet}
\int_{\Omega}\cu u \cdot \cu \varphi\, dx + \theta \, \int_{\Omega }\di u\,  \di \varphi \, dx = A  \int_{\Omega}u \cdot \varphi \, dx\,  \,, \forall\ \varphi \in (H^1_0(\Omega))^3,
\end{equation}
in the unknowns $u\in (H^1_0(\Omega))^3$ (the Dirichlet eigenfunction) and $A\in \mathbb{R}$ (the Dirichlet eigenvalue).

\begin{thm}\label{existence1}
Let $\Omega$ be a bounded, connected open set in $\R^3$ with  $C^{1,1}$ boundary. Let $\alpha <A_1$ and $\theta >0$.  Then there exists $c_{\alpha , \theta }\geq 0$ such that for any $\eta\geq c_{\alpha , \theta}$ the quadratic form defined by the left-hand side of \eqref{weakcucu1} is coercive in $X_{\rm \scriptscriptstyle\tiny T}(\Omega)$ and
problem \eqref{weakcucu1}  has a unique solution $u\in X_{\rm \scriptscriptstyle\tiny T}(\Omega)$ for all $f \in TL^2(\Gamma )$. Moreover, for $\alpha\le 0 $ one can take $ c_{\alpha , \theta}=0$.  Finally,  if $\alpha\le 0$ and 
$\eta=0$ and if, in  addition, $f$ satisfies the condition $\dig\, f=0$ on $\Gamma$, then $\di u=0$ in $\Omega$.
\end{thm}

{\bf Proof.}  We first consider the case  $\alpha \le 0$.  By a straightforward application of the Riesz Theorem combined with Theorem~\ref{gaffney}, and by noting that
the right-hand side of equation  \eqref{weakcucu1} defines linear continuous operators belonging to the dual space of $X_{\rm \scriptscriptstyle\tiny T}(\Omega )$, it is clear that the
problem is uniquely solvable in $X_{\rm \scriptscriptstyle\tiny T}(\Omega)$ for any $\eta \geq 0$.

We now consider the case $\alpha >0$. Let $\beta\in (\alpha , A_1)$. 
First of all we note that there exists $M>0$ such that 
\begin{equation}\label{poin}
 \int_{\Omega}|\cu u|^2\, dx + \theta \, \int_{\Omega }|\di u|^2 \, dx  +M\int_{\Gamma}|u^2|d\sigma \geq \beta \int_{\Omega} |u|^2\, dx  , 
\end{equation}
for all $u\in X_{\rm \scriptscriptstyle\tiny T}( \Omega)$.  Indeed, assume for contradiction that 
 for any $n\in {\mathbb{N}}$ there exists $u_n\in X_{\rm \scriptscriptstyle\tiny T}( \Omega)$ such that 
 \begin{equation}\label{poin1}
 \int_{\Omega}|\cu u_n|^2\, dx + \theta \, \int_{\Omega }|\di u_n|^2 \, dx  +n\int_{\Gamma}|u^2_n|d\sigma \le \beta \int_{\Omega} |u_n|^2\, dx  , 
\end{equation}
 and normalize $u_n$ by setting $\int_{\Omega}|u_n|^2 \, dx=1$. By \eqref{poin1} the sequence $u_n$, $n\in {\mathbb{N}}$, is bounded in $X_{\rm \scriptscriptstyle\tiny T}( \Omega)$; 
 hence - possibly passing to a subsequence - there exists $u \in X_{\rm \scriptscriptstyle\tiny T}( \Omega)$ such that $u_n\to u$ weakly in $ X_{\rm \scriptscriptstyle\tiny T}( \Omega)$, $u_n\to u $
 strongly in $L^2(\Omega)$. In particular, $\int_{\Omega}|u|^2 \, dx=1$. Moreover, since the trace operator is compact,  using condition 
 \eqref{poin1} allows to conclude that the trace of $u$ is zero, hence $u\in (H^1_0(\Omega))^3$. By passing to the limit in inequality \eqref{poin1} and using  the weak lower semicontinuity of 
 norms, we conclude that 
 \begin{eqnarray}
 \lefteqn{
  \int_{\Omega}|\cu u|^2\, dx + \theta \, \int_{\Omega }|\di u|^2 \, dx  }\nonumber \\
& &   \le \liminf_{n\to \infty } \left( \int_{\Omega}|\cu u_n|^2\, dx + \theta \, \int_{\Omega }|\di u_n|^2 \, dx\right)\nonumber \\
& &  \le  \lim_{n\to \infty }   \beta \int_{\Omega} |u_n |^2\, dx =\beta 
 \end{eqnarray}
 which implies that the infimum in \eqref{ray} is not larger than $\beta$, hence it is strictly smaller than $A_1$, a contradiction. 
 
 Let $\epsilon \in (0,1-  \frac{\alpha}{\beta} )$. By using \eqref{poin} we deduce immediately that 
 \begin{eqnarray}\lefteqn{
 \int_{\Omega}|\cu u|^2\, dx + \theta \, \int_{\Omega }|\di u|^2 \, dx  - \alpha \int_{\Omega} |u|^2\, dx   +(1-\epsilon )M\int_{\Gamma}|u^2|d\sigma} \nonumber \\
 & & 
= \epsilon \left(\int_{\Omega}|\cu u|^2\, dx + \theta \, \int_{\Omega }|\di u|^2 \, dx\right)  +
(1-\epsilon )\left(\int_{\Omega}|\cu u|^2\, dx  \right. \nonumber  \\
& & \left. +\, \theta \, \int_{\Omega }|\di u|^2 \, dx   +M\int_{\Gamma}|u^2|d\sigma\right)  - \alpha \int_{\Omega} |u|^2\, dx  \nonumber  \\
& & \geq  \epsilon \left(\int_{\Omega}|\cu u|^2\, dx + \theta \, \int_{\Omega }|\di u|^2 \, dx\right)  +((1-\epsilon )\beta -\alpha   )  \int_{\Omega} |u|^2\, dx\,,
\end{eqnarray}
   for all $u\in X_{\rm \scriptscriptstyle\tiny T}( \Omega)$, which implies the coercivity in  $X_{\rm \scriptscriptstyle\tiny T}( \Omega)$ (hence in 
$(H^1(\Omega))^3$) of the quadratic form in the left-hand side of \eqref{weakcucu1}. The conclusion about the existence of a unique solution then follows choosing
$c_{\alpha, \theta }=(1-\epsilon )M$ and proceeding as in the case $\alpha <0$.

 Now, assume that $\dig\, f=0$ on $\Gamma$, $\alpha\le 0$ and $\eta=0$.
Consider a solution $\Phi \in H^2(\Omega)$ to the following problem
\begin{equation}
\left\{  \begin{array}{ll}
\Delta \Phi = \di u,& \ \ {\rm in }\ \Omega\, , \\
D_{\nu } \Phi =0, &  \ \ {\rm on }\ \Gamma\, .
\end{array}
\right.
\end{equation}
Let us first observe that the above problem has a solution $\Phi  \in H^1(\Omega)$, since $u\cdot \nu=0$ on $\Gamma$ implies that   $\int_{\Omega }\di u\, dx =0$; furthermore, $\Phi \in H^2(\Omega)$ because $u\in H^1(\Omega)$ and therefore $\di u\in L^2(\Omega)$.  Clearly, $\gr \Phi \in  X_{\rm \scriptscriptstyle\tiny T}(\Omega)$ hence we can use
$\overline{\gr \Phi }$ as a test function in \eqref{weakcucu1} and get
\begin{eqnarray}\label{existence10}\lefteqn{-\alpha \int_{\Omega}|\gr \Phi|^2 dx+
\theta \int_{\Omega } | \di u  |^2\, dx  =\int_{\Omega }\cu u \cdot \cu \, \overline{\gr \Phi}\, dx  }\nonumber \\
& & -\, \alpha \int_{\Omega } u \cdot \overline{\gr \Phi}  \, dx + \theta \int_{\Omega }  \di u\, \overline{\Delta \Phi}\, dx =-\int_{\Gamma} f \cdot  \overline{\gr \Phi}\, d\sigma  =0\, ,
\end{eqnarray}
where the last equality holds because $\dig\, f=0$ on $\Gamma$, in view of the fact that $f$ is tangential and  hence $ f\,  \gr \Phi = f\,    \grg \Phi $ on $\Gamma$.
It follows from \eqref{existence10} and the assumption $\alpha\le 0$  that $\di u=0$ in $\Omega$.\hfill $\Box$

\begin{rem}
We note that if $\alpha \le 0$, $\theta>0$,  $\eta =0$ and $\dig\, f=0$ on $\Gamma$, then problem \eqref{weakcucu1} can be formulated directly in the energy space $X_{\rm \scriptscriptstyle\tiny T}(\di\, 0, \Omega )$. Therefore the weak formulation can be stated as: find
 $u\in X_{\rm \scriptscriptstyle\tiny T}(\di\, 0, \Omega )$ such that
\begin{equation}\label{weakcucu2}
\int_{\Omega}\cu u \cdot \cu \varphi\, dx  -  \alpha \int_{\Omega}u \cdot \varphi\, dx =-\int_{\Gamma} f \cdot \varphi \, d\sigma\, ,
\end{equation}
for all $\varphi \in X_{\rm \scriptscriptstyle\tiny T}(\di\, 0,\Omega )$. 

In order to see that,
 under these assumptions,  problems
\eqref{weakcucu1} and \eqref{weakcucu2} are equivalent, we can argue as follows. 
Consider a fixed $f \in TL^2(\Gamma )$  with $\dig\, f=0$ on $\Gamma$. By Theorem~\ref{existence1} the solution $u$ of problem
\eqref{weakcucu1} is also a solution of problem \eqref{weakcucu2}. To prove the reverse statement, consider a solution $u$ of \eqref{weakcucu2} and a test function
$\varphi \in X_{\rm \scriptscriptstyle\tiny T}(\Omega)$.
Along the same lines of the proof of Theorem~\ref{existence1},
consider a solution $\Psi \in H^2(\Omega)$ to the following problem
\begin{equation}
\left\{  \begin{array}{ll}
\Delta \Psi = \di \varphi,& \ \ {\rm in }\ \Omega\, , \\
D_{\nu } \Psi =0, &  \ \ {\rm on }\ \Gamma\, .
\end{array}
\right.
\end{equation}
Then we can decompose $\varphi$ as
$
\varphi =\widetilde \varphi +\gr \Psi \, ,
$
where  $\widetilde \varphi :=\varphi -\gr \Psi $. Since $\di  \widetilde \varphi = 0$ and  also $\di  u = 0$, we have that
\begin{equation}\label{equivalence0}
\int_{\Omega}\cu u \cdot \cu \varphi\, dx - \alpha\int_{\Omega}u \cdot \varphi\, dx+  \theta \int_{\Omega }\di u\,  \di \varphi \, dx = -\alpha   \int_{\Omega}u \cdot \widetilde \varphi\, dx +
    \int_{\Omega}\cu u \cdot \cu \widetilde \varphi\, dx \, .
\end{equation}
Moreover, since  $\dig\, f=0$ on $\Gamma$ we have
\begin{equation}\label{equivalence1}
\int_{\Gamma} f \cdot \varphi \, d\sigma = \int_{\Gamma} f \cdot (\widetilde \varphi + \gr \Psi )\, d\sigma  = \int_{\Gamma} f \cdot \widetilde \varphi \, d\sigma\, .
\end{equation}
By \eqref{equivalence0} and \eqref{equivalence1} we have that
\begin{equation}\label{equivalence2}
\int_{\Omega}\cu u \cdot \cu \varphi\, dx - \alpha\int_{\Omega}u \cdot \varphi\, dx + \theta \int_{\Omega }\di u\,  \di \varphi \, dx=-\int_{\Gamma} f \cdot \varphi \, d\sigma\, ,
\end{equation}
hence $u$ is a solution of problem \eqref{weakcucu1}. 
\end{rem}

\subsection{Resolvent operators}
For any $\alpha<A_1 , \theta >0$ and $\eta \geq c_{\alpha , \theta}$ as in Theorem~\ref{existence1},  we consider the operator $\mathcal{L}_{\alpha , \theta}^{\eta }$ from $X_{\rm \scriptscriptstyle\tiny T}(\Omega )$ to its dual $(X_{\rm \scriptscriptstyle\tiny T}(\Omega ))'$
defined by the pairing
$$
\langle \mathcal{L}_{\alpha, \theta}^{\eta }(u), \varphi \rangle = \int_{\Omega}\cu u \cdot \cu \varphi\, dx -\alpha \int _{\Omega}u \cdot \varphi\, dx+   \theta \int_{\Omega }\di u\,  \di \varphi \, dx
 +\eta \int_{\Gamma}u \cdot \varphi \, d\sigma
\, ,
$$
for all $u, \varphi \in X_{\rm \scriptscriptstyle\tiny T}(\Omega)$.  

Next, we consider the operator $\mathcal{J}$ from $TL^2(\Gamma)$ to $(X_{\rm \scriptscriptstyle\tiny T}(\Omega ))'$
defined by the pairing
$$
\langle \mathcal{J}(f), \varphi \rangle = \int_{\Gamma} f \cdot \varphi\,  d\sigma\, ,
$$
for all $f\in TL^2(\Gamma)$, $\varphi \in X_{\rm \scriptscriptstyle\tiny T}(\Omega)$.   

Recall that $\pi_{\rm \scriptscriptstyle\tiny T}(u)$ is the tangential components trace of $u$, which coincides with the trace of $u$ on $\Gamma$ for any $u\in X_{\rm \scriptscriptstyle\tiny T}(\Omega)$.
By Theorem~\ref{existence1},  $\mathcal{L}_{\alpha, \theta}^{\eta}$ is invertible, hence we can introduce the operator ${\mathcal{A}}^{\scriptscriptstyle\Gamma}_{\eta}$ defined from  $TL^2(\Gamma)$ to itself, by
$$
{\mathcal{A}}^{\scriptscriptstyle\Gamma}_{\eta}:=-\pi_{\rm \scriptscriptstyle\tiny T}\circ \left( \mathcal{L}_{\alpha, \theta}^{\eta} \right)^{-1} \circ \mathcal{J}\, .
$$

\begin{thm}
The operator ${\mathcal{A}}^{\scriptscriptstyle\Gamma}_{\eta} $ is compact and self-adjoint in $TL^2(\Gamma)$.
\end{thm}

{\bf Proof.} The compactness of  ${\mathcal{A}}^{\scriptscriptstyle\Gamma}_{\eta} $  follows from the fact that  the classical trace operator from $(H^1(\Omega))^3$ to $(L^2(\Gamma))^3$ is compact and coincides  with  the operator $\pi_{\rm \scriptscriptstyle\tiny T}$ on $X_{\rm \scriptscriptstyle\tiny T}(\Omega)$.  Regarding self-adjointness, it suffices to note that for all $f,g\in TL^2(\Gamma)$ we
have
\begin{eqnarray}\lefteqn{
\left\langle \left( \mathcal{L}_{\alpha, \theta}^{\eta} \right)^{-1}\circ \mathcal{J}(f), g \right\rangle_{(L^2(\Gamma))^3} =   \left\langle \mathcal{J}(\bar g), \left( \mathcal{L}_{\alpha, \theta}^{\eta} \right)^{-1} \circ \mathcal{J}(f) \right\rangle}\nonumber \\
& &  = \left\langle \mathcal{L}_{\alpha, \theta}^{\eta} \left(  \left( \mathcal{L}_{\alpha, \theta}^{\eta} \right)^{-1}   \circ \mathcal{J}(\bar g) \right),  \left( \mathcal{L}_{\alpha, \theta}^{\eta} \right)^{-1}  \circ \mathcal{J}(f) \right\rangle\nonumber \\
& & =\left\langle \mathcal{L}_{\alpha, \theta}^{\eta} \left( \left( \mathcal{L}_{\alpha, \theta}^{\eta} \right)^{-1} \circ \mathcal{J}(f) \right),  \left( \mathcal{L}_{\alpha, \theta}^{\eta} \right)^{-1}   \circ \mathcal{J}(\bar g) \right\rangle \nonumber \\
& & = \left\langle \mathcal{J}(f), \left( \mathcal{L}_{\alpha, \theta}^{\eta} \right)^{-1}  \circ \mathcal{J}(\bar g) \right\rangle = \left\langle f, \left( \mathcal{L}_{\alpha, \theta}^{\eta} \right)^{-1}   \circ \mathcal{J}(g) \right\rangle_{(L^2(\Gamma))^3} \,.
\end{eqnarray}
\hfill $\Box$\\

For our purposes, it is also convenient to consider the operator $  {\mathcal{A}}^{\scriptscriptstyle\Omega}_{\eta}$ from $X_{\rm \scriptscriptstyle\tiny T}(\Omega )$ to itself defined by
\begin{equation}
 {\mathcal{A}}^{\scriptscriptstyle\Omega}_{\eta}(u) = -\left( \mathcal{L}_{\alpha, \theta}^{\eta} \right)^{-1}  \circ \mathcal{J} \circ \pi_{\rm \scriptscriptstyle\tiny T}\, .
\end{equation}
It is evident that
\begin{equation}\label{link}
{\mathcal{A}}^{\scriptscriptstyle\Gamma}_{\eta}\circ \pi_{\rm \scriptscriptstyle\tiny T}=\pi_{\rm \scriptscriptstyle\tiny T}\circ  {\mathcal{A}}^{\scriptscriptstyle\Omega}_{\eta} .
\end{equation}

For $\eta\geq 0$,  it is  convenient to define the following sesquilinear form
\begin{equation}\label{sesq}
\langle u, v\rangle_{\alpha , \theta}^{\eta }:= \int_{\Omega}\cu u \cdot \cu \bar v \, dx -\alpha \int_{\Omega }u \cdot \bar v \, dx   + \theta \int_{\Omega }\di u\,  \di \bar v \, dx  +\eta\int_{\Gamma}u \cdot \bar v \, d\sigma   \, ,
\end{equation}
for all $u, v\in   X_{\rm \scriptscriptstyle\tiny T}( \Omega)$      and to observe that, if  $\eta\geq c_{\alpha , \theta}$,  it defines a scalar product in $X_{\rm \scriptscriptstyle\tiny T}(\Omega )$ by Theorem~\ref{existence1}, see also Theorem~\ref{gaffney}. In view of this,  problem \eqref{weakcucu1} can be written as 
$$
\langle u, \varphi \rangle_{\alpha , \theta}^{\eta}= -\langle f , \varphi \rangle_{(L^2(\Gamma))^3}   \, , 
$$
for all $\varphi\in  X_{\rm \scriptscriptstyle\tiny T}(\Omega )$, where $\langle \cdot , \cdot \rangle_{(L^2(\Gamma))^3}$ denotes the standard scalar product in   $(L^2(\Gamma))^3$     defined by $ \langle    f, \varphi \rangle_{(L^2(\Gamma))^3} :=     \int_{\Gamma}u \cdot \bar \varphi \, d\sigma $.

Note that  $\langle \cdot , \cdot\rangle_{\alpha , \theta}^0$ is the sesquilinear form appearing in the left-hand side of  \eqref{weakcucu1e}.

Then we have the following result, the proof of which is similar to the one of the previous theorem.

\begin{thm}
The operator $ {\mathcal{A}}^{\scriptscriptstyle\Omega}_{\eta} $ is compact and self-adjoint with respect to \eqref{sesq}.
\end{thm}

\subsection{The eigenvalue problem}
\label{eigensubsec}
In this subsection we consider the eigenvalue problem
 \eqref{classiccucu1e}. Recall that the weak formulation of  \eqref{classiccucu1e} is given in
\eqref{weakcucu1e}.
 It turns out that this eigenvalue problem can be recast as an eigenvalue problem for the operator $ {\mathcal{A}}^{\scriptscriptstyle\Gamma}_{\eta}$,  or for the  operator $  {\mathcal{A}}^{\scriptscriptstyle\Omega}_{\eta}$.
 Since these operators are compact and self-adjoint, their spectra can be easily described. In particular, we have the following result. 
   \begin{thm}\label{compactspectrum}  Let $\alpha <A_1$ and $\theta >0$. The spectrum  of the operator $ {\mathcal{A}}^{\scriptscriptstyle\Omega}_{\eta}$  can be represented as $\{ 0\} \cup \{ \gamma_n:\,  n\in\N\}  $
  where $\gamma_n,\ n\in \N $, are negative
   eigenvalues of finite multiplicity, $\gamma_n\to 0$ as $n\to \infty $, and $0$ is an eigenvalue of infinite multiplicity with eigenspace given by $(H^1_0(\Omega ))^3$. Moreover, the point spectrum of the operator ${\mathcal{A}}^{\scriptscriptstyle\Gamma}_{\eta}$  is given by $ \{\gamma_n:\, n\in\N\}$.  Furthermore, if $ {\mathcal{A}}^{\scriptscriptstyle\Omega}_{\eta}u=\gamma_n u$ for some $u\in X_{\rm \scriptscriptstyle\tiny T}(\Omega)$, then
  ${\mathcal{A}}^{\scriptscriptstyle\Gamma}_{\eta}\pi_{\rm \scriptscriptstyle\tiny T}u=\gamma_n \pi_{\rm \scriptscriptstyle\tiny T}u$.
  \end{thm}

{\bf Proof. } It is easy to see that the operators $ {\mathcal{A}}^{\scriptscriptstyle\Omega}_{\eta}$, ${\mathcal{A}}^{\scriptscriptstyle\Gamma}_{\eta}$ have the same non-zero eigenvalues and that zero is an eigenvalue of infinite multiplicity for the operator $ {\mathcal{A}}^{\scriptscriptstyle\Omega}_{\eta}$ with eigenspace equal to
$(H^1_0(\Omega ))^3$. The rest of the proof follows by the Hilbert-Schmidt Theorem applied to the compact and self-adjoint operator  ${\mathcal{A}}^{\scriptscriptstyle\Gamma}_{\eta}$ and by implementing \eqref{link}. \hfill $\Box$\\

Next, by Theorem~\ref{compactspectrum} and using the Min-Max Principle for the compact, self-adjoint operator $ {\mathcal{A}}^{\scriptscriptstyle\Omega}_{\eta}$,  we obtain the following result.

  \begin{thm}\label{mainspectrum}   Let $\alpha <A_1$ and $\theta >0$.
  The  eigenvalues of problem \eqref{classiccucu1e} form a sequence $\lambda_n,\ n \in \N$, in $ \R $, given by $\lambda_n=   \gamma_n^{-1}+\eta $, for all $n\in\N $  and the  eigenfunctions coincide with those of the operator  $ {\mathcal{A}}^{\scriptscriptstyle\Omega}_{\eta}$ associated with $\gamma_n$.
  Moreover,  $\lambda_n\to -\infty $, as $n\to \infty$, and can be represented as
  \begin{equation}\label{minmax}
  \lambda_n=  - \min_{ \substack{ V\subset X_{\rm \scriptscriptstyle\tiny T}(\Omega )  \\ {\rm dim }V=n }  }\  \,  \max _{u\in V\setminus (H^1_0(\Omega ))^3}
  \frac{  \int_{\Omega} \left( |\cu u|^2 -\alpha  |u|^2 + \theta  |\di u |^2 \right) dx }{\int_{\Gamma} |\pi_{\rm \scriptscriptstyle\tiny T}u|^2\, dx}\, ,
  \end{equation}
  where, as usual, each eigenvalue is repeated as many times as its multiplicity.
  \end{thm}

{\bf Proof.} It is enough to observe   that  $u\in X_{\rm \scriptscriptstyle\tiny T}(\Omega ) $ and $\lambda \in {\mathbb{R}}$
 satisfy equation   \eqref{weakcucu1e} if and only if  
 $$
 \langle u, \varphi\rangle_{\alpha , \theta}^{\eta}=-(\lambda -\eta)   \langle u, \varphi \rangle_{(L^2(\Gamma))^3}
 $$
 for all $\varphi \in X_{\rm \scriptscriptstyle\tiny T}(\Omega ) $,
 and this holds if and only if
$\lambda -\eta<0$ and
$
 {\mathcal{A}}^{\scriptscriptstyle\Omega}_{\eta}u=\gamma u
$,
where $\gamma =  (\lambda -\eta)^{-1}$. By Theorem~\ref{compactspectrum} we deduce the existence of the sequence of eigenvalues  $\lambda_n$, $n\in \mathbb{N}$,
by the sequence $\gamma_n$, $n\in {\mathbb{N}}$.  Since any eigenvalue $\lambda_n $ is bounded above by $\eta $ and $\gamma_n\to 0$, we deduce that
$\lambda_n\to -\infty $ as $n\to \infty$.  Formula \eqref{minmax} follows by the Min-Max Principle applied to the operator $- {\mathcal{A}}^{\scriptscriptstyle\Omega}_{\eta}$
 in which case one eventually
obtains the min-max characterization for $-(\lambda_n-\eta )$ which yields \eqref{minmax}. \hfill $\Box$\\

We note that the eigenvalue problem for the operator  ${\mathcal{A}}^{\scriptscriptstyle\Omega}_{\eta}$, can be written in the form
$$ \int_{\Omega}\cu u \cdot \cu \varphi\, dx -\alpha \int _{\Omega}u \cdot \varphi\, dx+   \theta \int_{\Omega }\di u\,  \di \varphi \, dx
 +\eta \int_{\Gamma}u \cdot \varphi \, d\sigma  = -\lambda  \int_{\Gamma} u \cdot \varphi \, d\sigma , $$
for all $\varphi\in  X_{\rm \scriptscriptstyle\tiny T}(\Omega )$,  in the unknowns $u \in  X_{\rm \scriptscriptstyle\tiny T}(\Omega )$ (the eigenvector) and $\lambda$ (the eigenvalue).  
It follows by the previous results that the space $X_{\rm \scriptscriptstyle\tiny T}(\Omega)$ can be decomposed as an orthogonal sum with respect to the scalar product \eqref{sesq}, namely
$$
X_{\rm \scriptscriptstyle\tiny T}(\Omega )=  {\rm Ker}  {\mathcal{A}}_{\eta}  ^{\scriptscriptstyle\Omega} \oplus   \left( {\rm Ker}  {\mathcal{A}}_{\eta}  ^{\scriptscriptstyle\Omega} \right)^{\perp } = (H^1_0(\Omega ))^3\oplus   \left( {\rm Ker}  {\mathcal{A}}_{\eta} ^{\scriptscriptstyle\Omega} \right)^{\perp } \, .
$$
We note that
$u\in  \left( {\rm Ker}  {\mathcal{A}}_{\eta} ^{\scriptscriptstyle\Omega} \right)^{\perp} $
if and only if

\begin{equation}\label{armoniche}
\int_{\Omega}\cu u \cdot \cu \varphi\, dx -\alpha \int_{\Omega}u \cdot \varphi \, dx+ \theta \int_{\Omega }\di u\,  \di \varphi \, dx= 0,
\end{equation}
 for all $  \varphi \in   (H^1_0(\Omega ))^3$ or, equivalently,  for all $  \varphi \in   (C^{\infty}_c(\Omega ))^3$. 
 
 Thus,  $u\in \left( {\rm Ker}  {\mathcal{A}}_{\eta} ^{\scriptscriptstyle\Omega} \right)^{\perp} $
 if and only if $u$ is a weak solution in $(H^1(\Omega))^3$ of the problem

\begin{equation}\label{weakharmonic}
\left\{
\begin{array}{ll}
\cu\, \cu u -\alpha u -\theta\, \gr\,  \di u=0,& \ \ {\rm in}\ \Omega ,\\
\nu \cdot u = 0 ,& \ \ {\rm on}\ \Gamma .
\end{array}
\right.
\end{equation}

 By setting\footnote{These functions are the analogues of the   harmonic functions considered in \cite{Auch2006}.  } $${\mathcal {H}}(\Omega):=  \left( {\rm Ker}  {\mathcal{A}} _{\eta}^{\scriptscriptstyle\Omega} \right)^{\perp}\,,$$ we can write
\begin{equation}
\label{dec}
X_{\rm \scriptscriptstyle\tiny T}(\Omega )=  (H^1_0(\Omega ))^3 \oplus   {\mathcal {H}}(\Omega).
\end{equation}

 By these observations, we deduce the validity of the following 

\begin{corol}   The operator  ${\mathcal{A}}^{\scriptscriptstyle\Gamma}_{\eta}$  has  no kernel, that is $ {\rm Ker}  {\mathcal{A}}_{\eta} ^{\scriptscriptstyle\Gamma} =\{0\}$.
\end{corol}
{\bf Proof.}   Assume that $f\in  {\rm Ker}  {\mathcal{A}}_{\eta} ^{\scriptscriptstyle\Gamma}$, that is 
$
-\pi_{\rm \scriptscriptstyle\tiny T}\circ \left( \mathcal{L}_{\alpha, \theta}^{\eta} \right)^{-1} \circ \mathcal{J}(f)=0
$.  We set $$u= \left( \mathcal{L}_{\alpha, \theta}^{\eta} \right)^{-1} \circ \mathcal{J}(f)$$ and we observe that, in particular, $u$ satisfies equation \eqref{armoniche} for all 
$  \varphi \in   (H^1_0(\Omega ))^3$ (note that in this step of the proof we have used only test functions vanishing at the boundary). Thus  $u\in {\mathcal {H}}(\Omega)$. Since by assumption $\pi_{\rm \scriptscriptstyle\tiny T}u=0$, we have that $u\in  (H^1_0(\Omega ))^3$. Thus, by \eqref{dec} we have $u=0$.  Going back to the equation satisfied by $u$, we obtain that 
$ \int_{\Gamma}f \cdot \varphi \, d\sigma=0$   for all $\varphi \in   X_{\rm \scriptscriptstyle\tiny T}(\Omega )  $, and this implies that $f=0$ since $f$ is a tangential field (note that in this second step, we have used all test functions  $\varphi \in   X_{\rm \scriptscriptstyle\tiny T}(\Omega )$ in the weak formulation of the equation). \hfill $\Box$ 

\subsection{Interior Calder\'{o}n operator}
\label{calderonsec}

In this section we introduce  a Calder\'{o}n operator associated with the interior problem \eqref{weakcucu1}. In order to identify the appropriate condition under which
our Calder\'{o}n operator is well-defined we need the following result.

In the sequel we denote by  $\Sigma =\{\lambda_n:\, n\in {\mathbb{N}} \}$ the set of Steklov eigenvalues of problem \eqref{weakcucu1e}.

\begin{thm}\label{existenceauto}
Let $\alpha <A_1$ and $\theta >0$. Then the problem
\begin{equation}\label{classiccucu1+eig}
\left\{
\begin{array}{ll}
\cu\, \cu u -\alpha u -\theta\, \gr\,  \di u=0,& \ \ {\rm in}\ \Omega ,\\
\nu \times \cu u  =\lambda \, u+ f ,& \ \ {\rm on}\ \Gamma
\end{array}
\right.
\end{equation}
is uniquely solvable  in $ X_{\rm \scriptscriptstyle\tiny T}(\Omega)$ for all $f \in TL^2(\Gamma )$ if and only if  $\lambda \notin \Sigma$.
\end{thm}
{\bf Proof. } Note  that problem \eqref{classiccucu1+eig} can be written in the weak form as
\begin{equation}
\label{fred}
\mathcal{L}^{\eta}_{\alpha, \theta }(u) = (\eta-\lambda ) \mathcal{J}(u) - \mathcal{J}(f)\, .
\end{equation}
Inverting the operator $\mathcal{L}^{\eta}_{\alpha, \theta }$,  problem \eqref{fred} turns out to be  equivalent to
\begin{equation}\label{fred1}
u=(\lambda -\eta){\mathcal{A}}_{\eta}^{ \scriptscriptstyle\Omega }(u)- \left( \mathcal{L}_{\alpha, \theta}^{\eta} \right)^{-1} \mathcal{J}(f)\,.
\end{equation}
Since ${\mathcal{A}}_{\eta}^{\scriptscriptstyle\Omega }$ is a compact operator, it follows by the Fredholm Alternative that problem \eqref{fred1} is uniquely solvable
if and only if $1$ is not an eigenvalue of the operator $(\lambda -\eta){\mathcal{A}}_{\eta}^{\scriptscriptstyle\Omega }$, and this exactly means  that
$\lambda\notin \Sigma$. \hfill $\Box$\\

Then   we can give the following definition.

\begin{defn}\label{calderondef} Assume that $\alpha  <A_1$ and $  \theta>0$ are such that $0\notin \Sigma$.  The interior Calder\'{o}n operator
$\cau$ is the operator defined from $TL^2(\Gamma)$ to itself, mapping any $f\in  TL^2(\Gamma)$
to $\cau (f) := \nu \times u $ on $\Gamma$, where $u\in X_{\rm \scriptscriptstyle\tiny T}(\Omega )$ is the solution of   \eqref{classiccucu1} (i.e., of  \eqref{weakcucu1} with $\eta =0$).
\end{defn}

We conclude this section by discussing the condition $0\notin \Sigma$. To do so, we consider two auxiliary eigenvalue problems. 
The first is the classical eigenvalue problem for the Neumann Laplacian
\begin{equation}\label{neulap}
\left\{ 
\begin{array}{ll}
-\Delta \phi =\lambda \phi,&\ \ {\rm in }\, \Omega,\vspace{1mm}\\
D_{\nu}\phi =0,&\ \ {\rm on }\, \partial\Omega\, ,
\end{array}
\right.
\end{equation}
for $\phi \in H^1(\Omega)$, 
which is well-known to admit a divergent sequence $\lambda_n^{{\scriptscriptstyle\mathcal{N}}}$, $n\in \N$, of non-negative eigenvalues of finite multiplicity, with $\lambda_1^{{\scriptscriptstyle\mathcal{N}}}=0$. 
The second is the eigenvalue problem for the $\cu\, \cu$  operator with ``magnetic'' boundary conditions
\begin{equation}\label{magbc}
\left\{
\begin{array}{ll}
\cu\, \cu u =\lambda u,& \ \ {\rm in}\ \Omega ,\\
\nu \times \cu u =0 ,& \ \ {\rm on}\ \Gamma ,
\end{array}
\right.
\end{equation}
for $u\in X_{\rm \scriptscriptstyle\tiny T}(\di \, 0, \Omega)$, 
which also   admits a divergent sequence $\lambda_n^{{\scriptscriptstyle\mathcal{M}}}$, $n\in \N$, of non-negative eigenvalues of finite multiplicity, with $\lambda_1^{{\scriptscriptstyle\mathcal{M}}}=0$, see e.g., \cite{zhang} for the relation between problem \eqref{magbc} and the eigenvalue problem for  Maxwell's system. 

Consider now problem \eqref{classiccucu1e} with $\lambda =0$, namely

\begin{equation}\label{classiccucu1ezero}
\left\{
\begin{array}{ll}
\cu\, \cu u -\alpha u -\theta\, \gr\,  \di u=0,& \ \ {\rm in}\ \Omega ,\\
\nu \times \cu u =0 ,& \ \ {\rm on}\ \Gamma .
\end{array}
\right.
\end{equation}
Then we can prove the following theorem where  we do not put any  {\it a priori} restrictions on $\alpha$ and $\theta$. For the values of $\alpha$ and $\theta$ for which $\Sigma $ is well-defined, the following statement gives a necessary and sufficient condition for the validity of the hypothesis $0\notin \Sigma $.  Note that the following theorem could be considered as the ``magnetic'' version of \cite[Theorem~1.1]{coda}. 

\begin{thm}  Problem   \eqref{classiccucu1ezero} has  a non-trivial solution $u\in X_{\rm \scriptscriptstyle\tiny T}(\Omega)$   if and only if $\alpha \in \{\theta \lambda_n^{{\scriptscriptstyle\mathcal{N}}}:\ n\in \N \} \cup  \{ \lambda_n^{{\scriptscriptstyle\mathcal{M}}}:\ n\in \N  \}$. 
\end{thm}

{\bf Proof.}  Since the case $\theta =0$ is trivial, we assume that $\theta \ne 0$. 

Clearly, for divergence free fields, \eqref{classiccucu1ezero} reads
\begin{equation}\label{classiccucu1ezerodiv}
\left\{
\begin{array}{ll}
\cu\, \cu u =\alpha  u,& \ \ {\rm in}\ \Omega ,\\
\nu \times \cu u =0 ,& \ \ {\rm on}\ \Gamma ,
\end{array}
\right.
\end{equation}
and  if $\alpha =\lambda_j^{{\scriptscriptstyle\mathcal{M}}}$ for some $j\in \N$, then the eigenfunctions corresponding to this $\lambda_j^{{\scriptscriptstyle\mathcal{M}}}$ solve problem 
\eqref{classiccucu1ezero}. 

Assume now  that  $\lambda=\frac{\alpha}{  \theta} $ is  an eigenvalue of problem \eqref{neulap}. Thus there exists 
a non-trivial solution $\phi \in H^1(\Omega)$  to  \eqref{neulap}.  In particular,  $\phi \in H^2(\Omega)$ by standard regularity theory. 
Then, it is readily seen that the function $u=\gr\, \phi $ is a non-trivial solution to \eqref{classiccucu1ezero}  if $\phi$ is not constant. 

 Thus, we have proved that if 
 $\alpha \in \{\theta \lambda_n^{{\scriptscriptstyle\mathcal{N}}}:\ n\in \N \} \cup  \{ \lambda_n^{{\scriptscriptstyle\mathcal{M}}}:\ n\in \N  \}$ then 
 there exists a non-trivial solution $u\in X_{\rm \scriptscriptstyle\tiny T}(\Omega)$  to \eqref{classiccucu1ezero}.  
 
 We now prove the converse statement.  Assume that 
  $u\in X_{\rm \scriptscriptstyle\tiny T}(\Omega)$ is a non-trivial solution   to \eqref{classiccucu1ezero} and assume that $\alpha \notin \{\theta \lambda_n^{{\scriptscriptstyle\mathcal{N}}}:\ n\in \N \}$. 
  Looking at the weak formulation of problem \eqref{classiccucu1ezero}, namely problem \eqref{weakcucu1e} with $\lambda =0$, and setting $\varphi =\gr \psi $ with $\psi \in H^2(\Omega)$ and $D_{\nu}\psi =0$, we get
  \begin{eqnarray}\label{calcvar}
  0= -\alpha\int_{\Omega}u \cdot \gr \psi \, dx +\theta\int_{\Omega}\di u\, \Delta \psi \, dx= \theta \int_{\Omega} \di u \left(\frac{\alpha}{\theta}\psi +\Delta \psi   \right)dx\, .
  \end{eqnarray}
  Since $\frac{\alpha }{ \theta }$  is assumed to be in the resolvent of the Neumann Laplacian,  the map $\frac{\alpha}{\theta}I +\Delta$ is a bijection
  between $\{\psi\in H^2(\Omega):\ D_{\nu}\psi =0 \}$ (the domain of the Neumann Laplacian for $\Omega$ of class $C^{1,1}$)  and  $L^2(\Omega)$.
 Thus, by the arbitrary choice of $\psi $ in \eqref{calcvar}, we deduce that $\di\,  u=0$, hence  
  $u$ is a solution of \eqref{classiccucu1ezerodiv}. In particular, $\alpha $ belongs to  $ \{ \lambda_n^{{\scriptscriptstyle\mathcal{M}}}:\ n\in \N  \}$. \hfill $\Box$

\section{Spectral representations}
\label{eigensec}

Let $\Omega$ be a bounded connected open set in $\R^3$ with boundary $\Gamma :=\partial \Omega \in C^{1,1}$.
Throughout this section we again assume that   $\alpha <A_1$ and $\theta >0$. 
 Recall that the space $X_{\rm \scriptscriptstyle\tiny T}(\Omega)$ can be considered as a Hilbert space with respect to the scalar product defined by \eqref{sesq}. Moreover, the space $ {\mathcal {H}}(\Omega)$ of solutions to problem \eqref{weakharmonic} is a closed subspace and it admits a Hilbert basis of Steklov eigenfunctions $u_n^{\scriptscriptstyle\Omega}$, $n\in \N$, which satisfy the equation
 \begin{equation}\label{weakcucu1en}
\int_{\Omega}\cu u_n^{\scriptscriptstyle\Omega} \cdot \cu \varphi\, dx -\alpha \int_{\Omega}u_n^{\scriptscriptstyle\Omega} \cdot \varphi \, dx+ \theta \int_{\Omega }\di u_n^{\scriptscriptstyle\Omega}\,  \di \varphi \, dx =  -  \lambda_n \int_{\Gamma} u_n^{\scriptscriptstyle\Omega} \cdot \varphi \, d\sigma\, ,
\end{equation}
for all $\varphi \in X_{\rm \scriptscriptstyle\tiny T}(\Omega)$.  Note that equation \eqref{weakcucu1en} can be equivalently written as 
 \begin{equation}\label{weakcucu1enbis}
 \langle u_n^{\scriptscriptstyle\Omega}, \varphi\rangle_{\alpha , \theta}^{\eta}=  
   -(\lambda_n-\eta)\langle u_n^{\scriptscriptstyle\Omega}, \varphi \rangle_{(L^2(\Gamma))^3}      \, ,
\end{equation}
for all $\varphi \in X_{\rm \scriptscriptstyle\tiny T}(\Omega)$.

In the sequel the eigenfunctions  $u_n^{\scriptscriptstyle\Omega}$, $n\in \N$, will be normalized with respect to the scalar product  \eqref{sesq}, namely we shall assume that 
$$
\left\langle u_n^{\scriptscriptstyle\Omega},  u_m^{\scriptscriptstyle\Omega} \right\rangle^{\eta}_{\alpha, \theta}= \delta_{nm}\,,
$$
where $\delta_{nm}$ is the
 Kronecker symbol. 
 
 Taking into account Theorem~\ref{compactspectrum} and the fact that the traces of the eigenfunctions $u_n^{\Omega }$  provide also a  basis of the space $TL^2(\Gamma)$, we set 
 \begin{equation}\label{norm}
u_n^{\scriptscriptstyle\Gamma}:= \sqrt{|\lambda_n-\eta |}\, \pi_{\rm \scriptscriptstyle\tiny T}u_n^{\scriptscriptstyle\Omega}
\end{equation}
and, in view of \eqref{weakcucu1enbis}, we observe that  $u_n^{\scriptscriptstyle\Gamma}$, $n\in \N$, is an orthonormal basis of $TL^2(\Gamma)$.

We proceed by proving the following theorem which  provides a spectral representation for the solutions of problem \eqref{classiccucu1} (i.e., of \eqref{weakcucu1} with $\eta =0$) for data $f\in TL^2(\Gamma) $ and  a corresponding spectral representation for the associated  Calder\'{o}n operator. 

Recall that, by Theorem~\ref{existenceauto}, if $0\notin \Sigma$ then problem  \eqref{classiccucu1} is uniquely solvable.  In Theorem~\ref{dualthm} we shall prove the same result for data $f\in TH^{-1/2}(\Gamma)$ and this will allow to  extend the Calder\'{o}n operator  and define it as an operator  from $TH^{-1/2}(\Gamma)$ to  $TH^{1/2}(\Gamma)$, as one would  expect.

\begin{thm}\label{representation} Assume that $0\notin \Sigma$. Let $f\in TL^2(\Gamma)$ be represented as
\begin{equation*}
f=\sum_{n=1}^{\infty }c_nu_n^{\scriptscriptstyle\Gamma}\, ,
\end{equation*}
where $(c_n)_{n\in \N}\in \ell^2$. Then the solution $u\in X_{\rm \scriptscriptstyle\tiny T}(\Omega)$ of problem  \eqref{classiccucu1} (i.e., of \eqref{weakcucu1} with $\eta =0$), is given by
\begin{equation}\label{representationu}
u= \sum_{n=1}^{\infty } \left(   \frac{\sqrt{|\lambda_n-\eta|}}{\lambda_n}\, c_n \right) u_n^{\scriptscriptstyle\Omega}\, .
\end{equation}
Moreover, the corresponding  interior Calder\'{o}n operator can be represented as
\begin{equation}\label{representationca}
\cau (f) = \nu \times \sum_{n=1}^{\infty }\frac{c_n}{\lambda_n}\, u_n^{\scriptscriptstyle\Gamma}\, .
\end{equation}
\end{thm}
{\bf Proof.}  Since $(c_n)_{n\in \N}\in \ell^2$ and $\lambda_n\to -\infty$, it is obvious that the series \eqref{representationu} converges in $X_{\rm \scriptscriptstyle\tiny T}(\Omega)$. By the continuity of
the form $\langle \cdot, \cdot\rangle_{\alpha , \theta}^0$ and of the trace operator, it follows that  it is  obvious that  
\begin{eqnarray}\lefteqn{
\langle u, \varphi \rangle_{\alpha , \theta}^0= \sum_{n=1}^{\infty } \left(   
\frac{\sqrt{|\lambda_n-\eta|}}{\lambda_n}\, c_n
\right)\! \left\langle u_n^{\scriptscriptstyle\Omega}, \varphi \right\rangle_{\alpha , \theta}^0 } \nonumber \\
& &=  -\sum_{n=1}^{\infty }   c_n\sqrt{|\lambda_n-\eta|}    
\langle u_n^{\scriptscriptstyle\Omega}, \varphi \rangle_{(L^2(\Gamma))^3 }  \nonumber  \\
& & 
=-\sum_{n=1}^{\infty} c_n     
\langle u_n^{\scriptscriptstyle\Gamma}, \varphi \rangle_{(L^2(\Gamma))^3}   
=       -
\langle f, \varphi \rangle_{(L^2(\Gamma))^3}
\, ,
\end{eqnarray}
 for all $\varphi\in X_{\rm \scriptscriptstyle\tiny T}(\Omega)$.  This means that $u$ is a solution of  problem \eqref{weakcucu1} with $\eta=0$.
Formula \eqref{representationca} follows immediately from  \eqref{representationu}.  \hfill $\Box$\\

In the spirit of \cite{Auch2006}, for $s>0$ we define the space $T{\mathcal{H}}^s(\Gamma)$ by  
\begin{equation}
\label{fracs}
T{\mathcal{H}}^s(\Gamma):=
\left\{f=\sum_{n=1}^{\infty }c_nu_n^{\scriptscriptstyle\Gamma}: \ \| f\|_{s,\Gamma}:=  \biggl(\sum_{n=1}^{\infty }|c_n|^2|\, \lambda_n-\eta|^{2s}    \biggr)^{1/2} < \infty    \right\}\,,
\end{equation}
endowed with the norm $\| \cdot \|_{s,\Gamma}$, and we define the space $T{\mathcal{H}}^{-s}(\Gamma)$ as the dual of $T{\mathcal{H}}^s(\Gamma)$.
Theorem~\ref{repdir}(i) shows that for $s=1/2$ and $s=-1/2$ this definition is equivalent to any of the classical definitions of the trace space: in other words
 $T{\mathcal{H}}^{1/2}(\Gamma)= T{{H}}^{1/2}(\Gamma)$ and $T{\mathcal{H}}^{-1/2}(\Gamma)= T{{H}}^{-1/2}(\Gamma)$, see Section~\ref{preliminaries}.

It is easy to see that the space $T{\mathcal{H}}^{-s}(\Gamma)$ can be identified with a space of sequences, namely
\begin{equation}
\left\{ F=(c_n)_{n\in \N}\in {\mathbb{C}}^{\N}:\     \| F\|_{-s,\Gamma}:= \biggl(\sum_{n=1}^{\infty }|c_n|^2\, |\lambda_n-\eta|^{-2s}   \biggr)^{1/2} < \infty     \right\} \,,
\end{equation}
with the understanding that the action of an element $F=(c_n )_{n\in \N }\in T{\mathcal{H}}^{-s}(\Gamma)$ on $f=\sum_{n=1}^{\infty }d_nu_n^{\scriptscriptstyle\Gamma}   \in  T{\mathcal{H}}^s(\Gamma )$ is given by the
pairing
\begin{equation}\label{action}
\langle F, f\rangle= \sum_{n=1}^{\infty } c_nd_n\,,
\end{equation}
which means that $\langle F, u_n^{\scriptscriptstyle\Gamma} \rangle = c_n$ for all $n\in \N$.

We note that any function $F\in TL^2(\Gamma)$ defines an element
of  $T{\mathcal{H}}^{-s}(\Gamma)$ by means of the formula
\[
\langle F, f\rangle = \int_{\Gamma } F \cdot f\,d\sigma\,\,,\,\,  \forall f\in T{\mathcal{H}}^{s}(\Gamma)\, .
\]
In this case,  in order to recover formula \eqref{action},  $F$ should be represented with res\-pect to the  basis  $\bar u_n^{\scriptscriptstyle\Gamma}:=  \overline {u_n^{\scriptscriptstyle\Gamma}}$, $n\in {\mathbb{N}}$.
 Indeed,  
if $F=\sum_{n=1}^{\infty }c_n\bar u_n^{\scriptscriptstyle\Gamma}$, then $\langle F, f\rangle = \int_{\Gamma } F \cdot f \,d\sigma = \sum_{n,m=1}^{\infty}c_nd_m\langle u_m^{\scriptscriptstyle\Gamma}, u_n^{\scriptscriptstyle\Gamma}\rangle_{(L^2(\Gamma))^3}= \sum_{n=1}^{\infty } c_nd_n$. Thus, we may think of $\bar u_n^{\scriptscriptstyle\Gamma}$, $n\in {\mathbb{N}}$, as a dual basis and represent formally any element
$F \in T{\mathcal H}^{-s}(\Gamma)$ as
$F=\sum_{n=1}^{\infty }c_n\bar u_n^{\scriptscriptstyle\Gamma}$. Note that since the coefficients of our operator are real,  it follows that if $u$ is an eigenfunction then also $\bar u$ is an eigenfunction  (see also the weak formulation \eqref{weakcucu1e}); hence one could choose a basis of real eigenfunctions in which case one should not worry about passing  from $u_n^{\scriptscriptstyle\Gamma}$ to its complex conjugate.)

As we have already mentioned in the introduction, the following theorem  allows in particular to characterize  the space $T{\mathcal{H}}^{1/2}(\Gamma)$ as the trace space of $X_{\rm \scriptscriptstyle\tiny T}(\Omega)$.  It  also  provides spectral representations for the solutions of the problem
\begin{equation}\label{weakcucudirbis}
\left\{
\begin{array}{ll}
\cu\, \cu u -\alpha u -\theta\, \gr\,  \di u=0,& \ \ {\rm in}\ \Omega ,\\
 u = f ,& \ \ {\rm on}\ \Gamma \, ,
\end{array}
\right.
\end{equation}
and of problem \eqref{weakcucudir}. We understand the solutions to problem \eqref{weakcucudirbis} as functions
$u\in {\mathcal{H}}(\Omega)$ satisfying the condition $\pi_{\rm \scriptscriptstyle\tiny T}u=u=f $. Note that by \eqref{dec} if $f=0$, the unique solution is $u=0$ since
$u$ would have to belong also  to the space $(H^1_0(\Omega))^3$. Thus, the solution to problem \eqref{weakcucudirbis} for all admissible data $f$ as below will be unique.

\begin{thm}\label{repdir}
The following statements hold
\begin{itemize}
\item[(i)]  The image of the trace operator $\pi_{\rm \scriptscriptstyle\tiny T}$ is given by
\begin{equation}
\pi_{\rm \scriptscriptstyle\tiny T}\left( X_{\rm \scriptscriptstyle\tiny T}(\Omega) \right)= T{\mathcal{H}}^{1/2}(\Gamma)\, ,
\end{equation}
hence $T{\mathcal{H}}^{1/2}(\Gamma)$ coincides with the usual Sobolev space  $T{H}^{1/2}(\Gamma)$.
\item[(ii)]
Let $f\in T{\mathcal{H}}^{1/2}(\Gamma)$ be represented as
\begin{equation}\label{repdir1}
f=\sum_{n=1}^{\infty }c_nu_n^{\scriptscriptstyle\Gamma}\,,
\end{equation}
where $\left( c_n\sqrt{|\lambda _n-\eta|} \right)_{n\in \N}\in \ell^2$. Then the solution $u\in X_{\rm \scriptscriptstyle\tiny T}(\Omega)$ of problem
\eqref{weakcucudirbis} is given by
\begin{equation}\label{repdir2}
u= \sum_{n=1}^{\infty}c_n\sqrt{|\lambda_n-\eta|}\, u_n^{\scriptscriptstyle\Omega}\, .
\end{equation}
\item[(iii)]
Let $f\in T{\mathcal{H}}^{1/2}(\Gamma)$ and let $f\times \nu$ be represented as
\begin{equation}\label{repdir3}
f\times \nu=\sum_{n=1}^{\infty } c_{n,\nu}\, u_n^{\scriptscriptstyle\Gamma}\,,
\end{equation}
where $\left( c_{n,\nu }\sqrt{|\lambda _n-\eta|} \right)_{n\in \N}\in \ell^2$. Then the solution $u\in X_{\rm \scriptscriptstyle\tiny T}(\Omega)$ of problem \eqref{weakcucudir} is given by
\begin{equation}\label{repdir4}
u= \sum_{n=1}^{\infty }c_{n,\nu }\, \sqrt{|\lambda_n-\eta|}\, u_n^{\scriptscriptstyle\Omega}\, .
\end{equation}
\end{itemize}
\end{thm}
{\bf Proof.} Let $u\in X_{\rm \scriptscriptstyle\tiny T}(\Omega)$. By \eqref{dec}  we have that $u=u_0+\widetilde u$, where $u_0\in (H^1_0(\Omega ))^3 $ and
$\widetilde u\in {\mathcal {H}}(\Omega)$. In particular, $\widetilde u$ can be written as
$\widetilde u= \sum_{n=1}^{\infty }  c_n u_n^{\scriptscriptstyle\Omega}$ with  $(c_n)_{n\in \N}\in \ell^2$.  Thus
$$\pi_{\rm \scriptscriptstyle\tiny T}u=\pi_{\rm \scriptscriptstyle\tiny T} \widetilde u = \sum_{n=1}^{\infty }  c_n  \pi_{\rm \scriptscriptstyle\tiny T} u_n^{\scriptscriptstyle\Omega}= \sum_{n=1}^{\infty} \frac{c_n}{\sqrt{|\lambda_n-\eta|}} u_n^{\scriptscriptstyle\Gamma}
\,,  $$
which clearly implies that $\pi_{\rm \scriptscriptstyle\tiny T}u \in  T{\mathcal{H}}^{1/2}(\Gamma) $. Thus, $\pi_{\rm \scriptscriptstyle\tiny T}\left( X_{\rm \scriptscriptstyle\tiny T}(\Omega)  \right)\subset T{\mathcal{H}}^{1/2}(\Gamma) $. 

In order to prove the
reverse inclusion, we consider a function $f\in T{\mathcal{H}}^{1/2}(\Gamma)$ represented as in
\eqref{repdir1} and we observe that the series in \eqref{repdir2} is convergent in $X_{\rm \scriptscriptstyle\tiny T}(\Omega)$. Thus, if $u$ is the function defined by \eqref{repdir2}, we have that
\begin{equation}\label{proofof2}
\pi_{\rm \scriptscriptstyle\tiny T}u =
\sum_{n=1}^{\infty} c_n\sqrt{|\lambda_n-\eta|}\,  \pi_{\rm \scriptscriptstyle\tiny T} u_n^{\scriptscriptstyle\Omega}
=\sum_{n=1}^{\infty }c_n u_n^{\scriptscriptstyle\Gamma}=f\,,
\end{equation}
which shows that $f\in \pi_{\rm \scriptscriptstyle\tiny T}\left( X_{\rm \scriptscriptstyle\tiny T}(\Omega) \right)$. Thus, $T{\mathcal{H}}^{1/2}(\Gamma)\subset \pi_{\rm \scriptscriptstyle\tiny T}\left( X_{\rm \scriptscriptstyle\tiny T}(\Omega) \right)$.

The equality \eqref{proofof2} gives also the proof of statement (ii), because, by definition,  the function $u$ in  \eqref{repdir2} belongs
to ${\mathcal{H}}(\Omega)$ hence it is a solution of \eqref{weakharmonic}. Moreover, since $u\cdot \nu =0$ on $\Gamma$, we have that
$\pi_{\rm \scriptscriptstyle\tiny T}u$ coincides with the trace of $u$, hence (with a slight abuse of notation) we have that $\pi_{\rm \scriptscriptstyle\tiny T}u=u$ on $\Gamma$. Thus, by
\eqref{proofof2}, the function  $u$ satisfies also the condition $u= f$ on $\Gamma$ required by problem \eqref{weakcucudirbis}, hence $u$ is a solution of \eqref{weakcucudirbis}.

The proof of statement (iii) follows directly by statement (ii) because the latter implies that the function $u$ defined by
\eqref{repdir4} is a solution of problem  \eqref{weakcucudirbis} with $f$ replaced by $f\times \nu$. Thus $u=f\times \nu$ on $\Gamma$. It follows that $\nu \times u=\nu\times (f\times \nu )= f $ as required in problem \eqref{weakcucudir}, hence the function $u$ defined by \eqref{repdir4} is a solution of problem  \eqref{weakcucudir}. \hfill $\Box$\\

Theorem~\ref{repdir} allows to consider equation \eqref{classiccucu1} also with a datum $f$ replaced by an element $F\in T{\mathcal{H}}^{-1/2}(\Gamma)$, in which case the formulation would read  as follows: find $u\in X_{\rm \scriptscriptstyle\tiny T}(\Omega)$ such that
\begin{equation}\label{dualeq}
\langle u, \varphi \rangle_{\alpha , \theta}^0=-\langle F, \bar\varphi \rangle\,,
\end{equation}
for all $\varphi \in  X_{\rm \scriptscriptstyle\tiny T}(\Omega) $. Indeed, the trace of $\varphi$ on $\Gamma$ belongs to $T{\mathcal{H}}^{1/2}(\Gamma)$ hence the right-hand side of equality \eqref{dualeq}
is well-defined. Note that the following theorem is stated with the use of an orthonormal basis of eigenfunctions which are not necessarily real. However, as we mentioned before, it is always  possible to select an orthonormal basis of real eigenfunctions. In accordance to the notation $\bar u_n^{\scriptscriptstyle\Gamma}=\overline{u_n^{\scriptscriptstyle\Gamma} }$ introduced before, we set
 $\bar u_n^{\scriptscriptstyle\Omega}:=\overline{u_n^{\scriptscriptstyle\Omega} }$.

\begin{thm}\label{dualthm} Assume that $0\notin \Sigma$.
Let $F\in T{\mathcal{H}}^{-1/2}(\Gamma)$ be represented as
\begin{equation}\label{dualthm1}
F=\sum_{n=1}^{\infty }c_n\bar u_n^{\scriptscriptstyle\Gamma}\,,
\end{equation} where $\left( c_n |\lambda_n-\eta|^{-1/2} \right)_{n\in \N}\in \ell^2$. Then the solution $u\in X_{\rm \scriptscriptstyle\tiny T}(\Omega)$ of problem \eqref{dualeq} is given by
\begin{equation}\label{dualthm2}
u= \sum_{n=1}^{\infty } \! \left(  \frac{\sqrt{|\lambda_n-\eta |}}{\lambda_n}\, c_n \right) \bar u_n^{\scriptscriptstyle\Omega}\, .
\end{equation}
\end{thm}
{\bf Proof.}  Note that  since $\left( c_n  |\lambda_n-\eta|^{-1/2} \right)_{n\in \N}\in \ell^2$,  the series \eqref{dualthm2} converges in $X_{\rm \scriptscriptstyle\tiny T}(\Omega)$.  Let $\bar \varphi =\sum_{n=1}^{\infty }d_nu_n^{\scriptscriptstyle\Gamma}\in T{\mathcal{H}}^{1/2}(\Gamma)$ with $\left( d_n  \sqrt{|\lambda_n-\eta|} \right)_{n\in N}\in \ell^2$.  By the continuity of
the bilinear form $\langle \cdot, \cdot\rangle_{\alpha , \theta}^0$ we have that 
\begin{eqnarray}\lefteqn{
\langle u, \varphi \rangle_{\alpha , \theta}^0= \sum_{n=1}^{\infty } \! \left(\!   \frac{\sqrt{|\lambda_n-\eta |}}{\lambda_n}\, c_n             \right)\! \left\langle \bar u_n^{\scriptscriptstyle\Omega}, \varphi \right\rangle_{\alpha , \theta}^0=\!  -\sum_{n=1}^{\infty }   c_n\sqrt{|\lambda_n-\eta|}      
\langle \bar u_n^{\scriptscriptstyle\Omega}, \varphi \rangle_{(L^2(\Gamma))^3}
 }
\nonumber \\
& &\qquad = - \sum_{n=1}^{\infty }  c_n\int_{\Gamma}\bar u_n^{\scriptscriptstyle\Gamma} \, \sum_{m=1}^{\infty}  d_m \, u_m^{\scriptscriptstyle\Gamma}  \,  d\sigma= -\sum_{n=1}^{\infty}  c_n \, d_n =- \langle F, \bar \varphi \rangle \,,
\end{eqnarray}
which  means that $u$ is a solution of  problem \eqref{dualeq}.  \hfill $\Box$

\begin{rem}\label{finalrem} By formula  \eqref{dualthm2}, it follows that the  interior Calder\'{o}n operator defined in Definition~\ref{calderondef} can be extended from $TL^2(\Gamma)$ to  $T{\mathcal{H}}^{-1/2}(\Gamma)$ by setting
\begin{equation}\label{dualthm3}
\cau (F) = \nu \times \sum_{n=1}^{\infty }\frac{c_n}{\lambda_n}\, \bar u_n^{\scriptscriptstyle\Gamma}\, ,
\end{equation}
for all $F  \in T{\mathcal{H}}^{-1/2}(\Gamma)$ represented as in \eqref{dualthm1},  with $\cau (F)$ being an element 
of   $T{\mathcal{H}}^{1/2}(\Gamma)$.
\end{rem}

\section{Appendix: the case $\alpha >A_1$}
\label{app}

Our approach allows to treat also the case $a> A_1$. First of all, one should note that by standard 
spectral theory, problem \eqref{dirichlet} has a divergent sequence of positive eigenvalues $A_n$, $n\in \mathbb{N}$
with finite multiplicity. Assume that $\alpha\in {\mathbb{R}}$ is such that 
$
A_n< \alpha <A_{n+1}
$ 
for some $n\in {\mathbb{N}}$.  
 
 Let $V_n$ be the subspace of $(H^1_0(\Omega))^3$ generated by all eigenfunctions associated with all eigenvalues $A_k$ with $k\le n$, and let 
 $$
 V_n^{\perp}=\{v\in X_{\rm \scriptscriptstyle\tiny T}(\Omega ) :\ \  \langle v, u\rangle_{\alpha, \theta  }^0=0,\ \forall \ u\in V_n  \}\, .
 $$ 
 
 Clearly, $V_n^{\perp}$ is a closed subspace of the $X_{\rm \scriptscriptstyle\tiny T}(\Omega )$.  Then we have the following result.

\begin{thm}\label{splitting}
Let $\Omega$ be a bounded, connected open set in $\R^3$ with  $C^{1,1}$ boundary. Assume that $
A_n< \alpha <A_{n+1}
$ 
for some $n\in {\mathbb{N}}$, and let $\theta >0$.  Then 
\begin{equation}\label{splitting1}
X_{\rm \scriptscriptstyle\tiny T}(\Omega )= V_n\oplus V_n^{\perp}\,, 
\end{equation} 
and  there exists $c_{\alpha , \theta }\geq 0$ such that for any $\eta \geq c_{\alpha , \theta}$ the quadratic form defined by the left-hand side of \eqref{weakcucu1} is coercive in $V_n^{\perp} $. 
\end{thm}
 
{\bf Proof. } We note that $\langle \cdot , \cdot\rangle_{\alpha, \theta  }^0$ is not necessarily a scalar product, hence the proof of \eqref{splitting1} requires 
some justification. 
Given $v\in  X_{\rm \scriptscriptstyle\tiny T}(\Omega )$, the weak problem 
$$
\langle u, \varphi \rangle_{\alpha , \theta  }^0=    \langle v, \varphi \rangle_{\alpha , \theta  }^0     ,\ \forall \ \varphi\in V_n \,, 
$$
has a solution $u\in V_n$. This  can be  proved by finding a critical point of the functional $u\mapsto \frac{1}{2}\langle u, u \rangle_{\alpha , \theta  }^0-    \langle v, u\rangle_{\alpha , \theta  }^0  $ in  the finite dimensional space $V_n$. 
To do so, note that for all $u\in V_n$ we have 
\begin{equation}\label{estvn} \int_{\Omega}|\cu u|^2\, dx + \theta \, \int_{\Omega }|\di u|^2 \, dx   \le A_n \int_{\Omega}|u|^2dx
\end{equation}
hence 
$
\langle u, u \rangle_{\alpha , \theta  }^0 \le -\varrho \int_{\Omega}|u|^2dx\,,
$
where $\varrho = \alpha-A_n>0$. Thus for $u\in V_n $, $u\ne 0$ we have
\begin{eqnarray}\lefteqn{
\frac{1}{2} \langle u, u \rangle_{\alpha , \theta  }^0-    \langle v, u\rangle_{\alpha , \theta  }^0}\nonumber \\
& & 
 \le  -\frac{\varrho}{2} \int_{\Omega}|u|^2dx-    \langle v, u\rangle_{\alpha , \theta  }^0
=
  -\frac{\varrho}{2} \| u\|_{L^2(\Omega)}^2-\| u\|_{L^2(\Omega)} \left\langle v, \frac{u}{ \| u\|_{L^2(\Omega)}  }\right\rangle_{\alpha , \theta  }^0  \nonumber \\
& &   \le 
  -\frac{\varrho}{2} \| u\|_{L^2(\Omega)}^2-\| u\|_{L^2(\Omega)}  \min_{  \substack{\| u\|_{L^2(\Omega)}  =1\\ u\in V_n     } }   \langle v,u\rangle_{\alpha , \theta  }^0\, ,
\end{eqnarray}
which implies that $\frac{1}{2} \langle u, u \rangle_{\alpha , \theta  }^0-    \langle v, u\rangle_{\alpha , \theta  }^0\to -\infty  $ as $\| u\|_{L^2(\Omega )} \to \infty $. Thus the map
 $u\mapsto \frac{1}{2}\langle u, u \rangle_{\alpha , \theta  }^0 -    \langle v, u\rangle_{\alpha , \theta  }^0  $ has actually a maximum in $V_n$, hence a critical point, as required. 

 Then $v= u+(v-u)$ with $v-u\in V_n^{\perp} $, hence $X_{\rm \scriptscriptstyle\tiny T}(\Omega )= V_n+V_n^{\perp}$. 
 
 On the other hand, the  relation  $V_n\cap V_n^{\perp}=\{0\}$ 
can be proved as follows. Assume by contradiction that  there exists $u\in V_n$ such that $ \langle u, \varphi \rangle_{\alpha }^0=0 $ for all $\varphi \in V_n$, then 
$\  \int_{\Omega}|\cu u|^2\, dx + \theta \, \int_{\Omega }|\di u|^2 \, dx = \alpha\int_{\Omega}|u|^2dx$    which opposes \eqref{estvn}.
Thus, \eqref{splitting1} holds. 

We claim that 
\begin{equation}\label{rayn}
\widetilde A_{n+1}:=\inf_{\substack{v\in  V_n^{\perp} \cap (H^1_0\Omega))^3  \\ v\ne 0 }  }\frac{  \int_{\Omega}|\cu v|^2\, dx + \theta \, \int_{\Omega }|\di v|^2 \, dx   }{\int_{\Omega} |v|^2dx}>\alpha .
\end{equation}
Indeed, assume for contradiction that $\widetilde A_{n+1}\le \alpha$. By standard arguments, we can find a minimizer $v \in V_n^{\perp}\cap H^1_0(\Omega)$ such that 
$ \int_{\Omega}|\cu v|^2\, dx + \theta \, \int_{\Omega }|\di v|^2 \, dx =\widetilde A_{n+1}    \int_{\Omega} |v|^2dx $ and we can assume that $\int_{\Omega}|v|^2dx=1$. Consider now the  $(n+1)$-dimensional space
$W_{n+1}=V_n+\langle v\rangle$. Given  an element $w=u+\xi v\in W_n$, with $\xi \in {\mathbb{C}}$, we easily see that 
\begin{eqnarray}\lefteqn{
  \int_{\Omega}|\cu w|^2\, dx + \theta \, \int_{\Omega }|\di w|^2 \, dx  =
\int_{\Omega}|\cu u|^2\, dx + \theta \, \int_{\Omega }|\di u|^2 \, dx   }\nonumber \\
& & 
+|\xi|^2\biggl(  \int_{\Omega}|\cu v|^2\, dx + \theta \, \int_{\Omega }|\di v|^2 \, dx       \biggr)   \nonumber  \\
& &+
2 {\rm Re}\, \bar\xi\biggl(\int_{\Omega}\cu u\, \cu \bar v \, dx   +    \theta \int_{\Omega }\di u\,  \di \bar v \, dx\biggr) \nonumber \\
& &= \int_{\Omega}|\cu u|^2\, dx + \theta \, \int_{\Omega }|\di u|^2 \, dx   +
|\xi|^2\left(  \int_{\Omega}|\cu v|^2\, dx + \theta \, \int_{\Omega }|\di v|^2 \, dx       \right) \nonumber \\
& & +
2 \alpha  {\rm Re}\, \bar\xi\int_{\Omega}u\bar vdx \le A_n \int_{\Omega}|u|^2dx+|\xi|^2 \widetilde A_{n+1}+2 \alpha  {\rm Re}\, \bar \xi  \int_{\Omega}u\bar vdx\nonumber  \\
& & 
\le  \alpha \left( \int_{\Omega}|u|^2dx+|\xi|^2 +2 {\rm Re}\, \bar\xi  \int_{\Omega}u\bar vdx  \right)= \alpha
 \int_{\Omega}|w|^2dx
\end{eqnarray}

This implies that 
$$
\sup_{\substack{w\in W_{n+1}\\w\ne 0}}  \frac{  \int_{\Omega}|\cu w|^2\, dx + \theta \, \int_{\Omega }|\di w|^2 \, dx}{ \int_{\Omega}|w|^2dx}\le \alpha
$$
and by the Min-Max Principle 
$$
A_{n+1}= \inf_{\substack{W\subset (H^1_0\Omega))^3\\ \dim W\le n+1}}  \sup _{\substack{w\in W\\w\ne 0}}  \frac{  \int_{\Omega}|\cu w|^2\, dx + \theta \, \int_{\Omega }|\di w|^2 \, dx}{ \int_{\Omega}|w|^2dx}\le \alpha
$$
which is clearly a contradiction since $\alpha <A_{n+1}$. The claim is proved. 

Inequality \eqref{rayn} allows to prove that 
 for a fixed $\beta \in (\alpha, \widetilde A_{n+1})$, there exists  $M>0$ such that 
\begin{equation}\label{poinenne}
 \int_{\Omega}|\cu u|^2\, dx + \theta \, \int_{\Omega }|\di u|^2 \, dx  +M\int_{\Gamma}|u^2|d\sigma \geq \beta \int_{\Omega} |u|^2\, dx  , 
\end{equation}
for all $u\in V_n^{\perp}$: the argument is the same one used in the proof of Theorem~\ref{existence1}, where the role of the space $X_{\rm \scriptscriptstyle\tiny T}(\Omega )$ is now played by  $V_n^{\perp}$. Again, by the same argument used in the proof of Theorem~\ref{existence1}, one can easily deduce by \eqref{poinenne} the coercivity of the quadratic form in the statement. \hfill $\Box$\vspace{12pt} \\
Now we
observe that, whenever $u\in V_n^{\perp} $, one 
can equivalently consider, in the weak formulation of problem \eqref{weakcucu1}, only test functions $\varphi \in V_n^{\perp}$;  indeed, adding to  $\varphi $ a function $\widetilde \varphi\in V_n$ leaves invariant both sides  of the equation. Therefore, all the analysis carried out for the various problems discussed in this paper, can still be performed without any essential modifications; it suffices to replace the energy space 
$X_{\rm \scriptscriptstyle\tiny T}(\Omega )$ by the smaller energy space $V_n^{\perp}$, thus excluding the Dirichlet eigenfunctions given by $V_n$.

\vspace{0.4cm}
\noindent
{\bf Acknowledgements}\\
 The authors are very thankful to an anonymous referee for the careful reading of the paper and many useful  remarks. They  are also very thankful  to Dr. Luigi Provenzano for stimulating discussions - on Steklov-type problems with negative potentials - which were inspirational for our observations contained in  Section~\ref{app}.
This work was performed while the second named author (IGS) visited the Department of Mathematics ``Tullio Levi-Civita'' of the University of Padova, Italy, in the framework of the 2018 Visiting Scientist Programme. IGS acknowledges that this work was made possible by NPRP grant \#[8-764-160] from Qatar National Research Fund (a member of Qatar Foundation).
 The first named author (PDL)  acknowledges partial  financial support from the research project BIRD191739/19  ``Sensitivity analysis of partial differential equations in the mathematical theory of electromagnetism'' of the University of Padova.  PDL  is also a member of the Gruppo Nazionale per l'Analisi Matematica, la Probabilit\`{a} e le loro Applicazioni (GNAMPA) of the Istituto Nazionale di Alta Matematica (INdAM), and acknowledges partial financial support from Progetto GNAMPA 2019  ``Analisi  spettrale  per  operatori  ellittici  con  condizioni  di  Steklov  o  parzialmente incernierate''.

$ $ \vspace{8pt}\\

\noindent Pier Domenico  Lamberti\\
Department of Mathematics ``Tullio Levi-Civita'', University of Padova,\\
Via Trieste 63, I-35121 Padova, Italy.\\
E-mail address: lamberti@math.unipd.it\\

\noindent Ioannis G. Stratis\\
Department of Mathematics, National and Kapodistrian University of Athens,\\
Panepistimiopolis, GR-15784 Zographou (Athens), Greece.\\
E-mail address: istratis@math.uoa.gr\\

\end{document}